\documentclass{amsart}

\usepackage{amssymb}
\usepackage{epsfig}
\usepackage{psfrag}
\newtheorem{theorem}{Theorem}[section]
\newtheorem{lemma}[theorem]{Lemma}
\newtheorem{proposition}[theorem]{Proposition}
\newtheorem{corollary}[theorem]{Corollary}

\theoremstyle{definition}
\newtheorem{definition}[theorem]{Definition}

{\bf }

\theoremstyle{remark}
\newtheorem{remark}[theorem]{Remark}
\numberwithin{equation}{section}
\newcommand{\abs}[1]{\lvert#1\rvert}

\newcommand{\ul}[1]{\underline{#1}}
\newcommand{\C}{\mathbb C}
\newcommand{\E}{\mathbb E}
\newcommand{\comment}[1]{} 

\begin{document}

\title{The planar algebra of group-type subfactors}

\author{Dietmar Bisch}
\address{Vanderbilt University, Department of Mathematics, SC 1326,
Nashville, TN 37240, USA}
%
\email{dietmar.bisch@vanderbilt.edu}
\email{paramita.das@vanderbilt.edu}
\email{shamindra.k.ghosh@vanderbilt.edu}

\thanks{The authors were supported by NSF under Grant No. DMS-0301173.}
%
\author{Paramita Das}
%
\author{Shamindra Kumar Ghosh}

%


\keywords{planar algebra, planar operad, subfactor, IRF model}

\begin{abstract}
If $G$ is a countable, discrete group generated by two finite
subgroups $H$ and $K$ and $P$ is a II$_1$ factor with an outer
G-action, one can construct the group-type subfactor $P^H \subset P \rtimes K$ 
introduced in \cite{BH}. This construction was used in \cite{BH} to obtain 
numerous examples of infinite depth subfactors whose standard invariant 
has exotic growth properties.
We compute the planar algebra (in the sense of Jones \cite{J2}) of this 
subfactor
and prove that any subfactor with an abstract planar algebra of 
"group type" arises from such a subfactor. The action of Jones' 
planar operad is determined explicitly.
\end{abstract}

\maketitle

\section{Introduction}
The technique of {\it composing subfactors} pioneered in \cite{BH} 
led to a zoo of exotic subfactors of the hyperfinite II$_1$ factor. 
In particular,
the first examples of irreducible, amenable subfactors that are not strongly
amenable (in the sense of \cite{P2}) were constructed in this way in
\cite{BH}.
The idea is simple: Let $H$ and $K$ be two finite groups with
an outer action on a II$_1$ factor $P$ (e.g. the hyperfinite II$_1$ 
factor) and consider the composition of the fixed-point subfactor 
$P^H \subset P$ with the crossed product subfactor 
$P \subset P \rtimes K$ to obtain, what we will call here, a
{\it group-type subfactor} $P^H \subset P \rtimes K$.
If $P$ is hyperfinite, analytical properties of this subfactor, such 
as amenability and property (T) in the sense of Popa (\cite{P2}, \cite{P4}), 
were proved to be equivalent to the corresponding properties 
(amenability, property (T) in the sense of Kazhdan) of the group
$G$ generated by $H$ and $K$ in the outer automorphism group of 
$P$ (\cite{BH}, \cite{BP}). Note that a group-type subfactor is obtained 
from the two fixed-point subfactors $P^H \subset P$ and $P^K \subset P$ 
by performing the basic construction of \cite{J1} to one of them. 
A group-type subfactor is therefore an invariant for the relative position 
of the two fixed-point subfactors. For more on this, see \cite{J5}.

The main invariant for a subfactor is the so-called {\it standard
invariant} (see for instance \cite{GHJ}, \cite{EK}, \cite{JS}). 
It is a very sophisticated mathematical object that can be
portrayed in a number of seemingly quite different ways. For example, 
it has descriptions as a certain category of bimodules ({\cite{Oc2}, 
see also \cite{Bi2}), as a lattice of algebras (\cite{P1}), or as a 
planar algebra (\cite{J2}). Jones' planar algebra technology has 
become a very efficient tool to capture and analyze the standard 
invariant of a subfactor. 

Composition of subfactors was the motivating idea that led to the
results in \cite{BJ1}. It was proved there that two standard
invariants without extra structure, i.e. consisting of just
the Temperley-Lieb algebras, can always be composed 
freely to form a new standard invariant -- namely the one generated 
by the Fuss-Catalan algebras of \cite{BJ1}. This concept of {\it free
composition} was then pushed much further in \cite{BJ2}, where
it is shown that any two planar algebras arising from subfactors
can be composed freely to form a new subfactor planar algebra.

The principal graphs of a subfactor encode the algebraic information 
contained in the standard invariant (\cite{GHJ}), and their structure 
determines the growth properties of the invariant. It was shown in
\cite{BH} that subfactors whose standard invariants have very
exotic growth properties exist by constructing concrete group-type
subfactors. The principal graphs of these subfactors were computed
there. In this paper we go a step further and give a concrete description 
of the standard invariant (or equivalently the planar algebra) of the
group-type subfactors. We concentrate on the case when the group
$G$, generated by $H$ and $K$ in the automorphism group
of $P$, has an outer action on $P$. Note that any group $G$ generated
by two finite subgroups $H$ and $K$ has such an action on the
hyperfinite II$_1$ factor (for instance by a Bernoulli shift action).
We find a description of the planar
algebra of these group-type subfactors that is reminiscent of
an IRF (interaction round a face) model in statistical mechanics.
The general case, where $G$ is generated by $H$ and $K$
in the {\it outer} automorphism group of $P$, and hence may or may not
lift to an action of $G$ on $P$, is much more elaborate and will be 
treated in a separate paper. 

Here is a more detailed outline of the sections of this paper.
We review in section 2 the basic notions of planar algebras. In
section 3, we define an abstract planar algebra $\mathcal{P}$ associated 
to a countable, discrete group $G$ and two of its finite subgroups
$H$ and $K$ which generate $G$. The vector spaces underlying the
planar algebra are spanned by alternating words in $H$ and $K$
that multiply to the identity element. The action of Jones'
planar operad is given explicity by a particular labelling
of planar tangles that can be viewed as an IRF-like model. 
It takes some work to show that the action of planar tangles
is well-defined and preserves composition. The latter is
achieved by showing that the action respects composition
with certain elementary tangles that generate any annular
tangle. We then analyze the filtered $*$-algebra structure
of $\mathcal{P}$ and determine the action of Jones projection and
conditional expectation tangles.

In section 4 we compute the basic construction and higher relative
commutants of a group-type subfactor $P^H \subset P \rtimes K$. 
We exhibit a nice basis which is used in section 5 to prove
that the abstract group-type planar algebra of section 3 is indeed 
isomorphic to the concrete one computed in section 4. Moreover, we prove
in section 5 that if the standard invariant of an arbitrary
subfactor is isomorphic to a group-type planar algebra, then
the subfactor is indeed of group-type. This is proved
using results on intermediate subfactors from \cite{Bi1}, \cite{BJ2},
\cite{BZ}.

\section{Planar algebra basics}\label{genplnalg}
In this section, we will give a very brief overview of 
planar algebras; the reader is encouraged to see \cite{J2} for details.
Let us first describe the key ingredients that constitute a 
{\em planar tangle}.
\begin{itemize}
\item There is an external disc, 
several (possibly none) internal discs and
a collection of disjoint smooth curves. 
\item To each disc - internal or external, we attach a 
non-negative integer. This integer will be referred to as the 
`color' of the disc.
If a disc has color $k > 0$, there will be $2k$ points on the 
boundary of the disc marked $1, 2, \cdots 2k$ counted clockwise, 
starting with a distinguished marked point, which is decorated with
`$*$'. 
A disc having color $0$ will have no marked points on its boundary.
\item Each of the curves is either closed, or joins a marked point on
the
boundary of a disc to another such point, meeting the boundary of the
disc transversally. Each marked point must be the endpoint of exactly 
one curve. 
\item The whole picture has to be planar, in the sense that there 
should be no crossing of curves or overlapping of discs. 
\item Finally, we will not distinguish between pictures obtained 
from one another by planar isotopy preserving the $*$'s. 
\end{itemize}
The data of such a picture will be termed as a {\em planar 
$k$-tangle}, where $k$ refers to the color of its external disc.
\begin{remark}
We can induce a black-and-white shading on the complement of the 
union of the internal discs and curves in the external disc by
specifying
that the region between the last and first marked point be left
unshaded. This leaves a scope for ambiguity in the case of
$0$-discs, thus, $0$-discs may be of two types depending on whether
the region surrounding their boundary is shaded or unshaded.    
\end{remark}
Given a planar $k$-tangle $T$ - one of whose internal discs have 
color $k_i$ -  and a $k_i$-tangle $S$, one can define the $k$-tangle 
$T \circ_i S$ by isotoping $S$ so that its boundary, together 
with the marked points and the `$*$' coincides with that of $D_i$ and
then
erasing the boundary of $D_i$. The collection of tangles - along with 
the composition defined thus - is called the {\em colored operad of
planar tangles}. 

A planar algebra is a collection of vector spaces $\{P_k\}_{k \geq 0}$
such that every $k$-tangle $T$ that has $b$ internal discs 
$D_1, D_2, \cdots D_b$ having colors $k_1, k_2, \cdots k_b$ 
respectively gives rise to a multilinear map 
$Z_T : P_{n_1} \times P_{n_2} \cdots \times P_{n_b} \rightarrow
P_{n_0}$. The collection of maps is required to be compatible 
with substitution of tangles and renumbering of internal discs.

\section{An abstract planar algebra}\label{plnalg}
In this section we will abstractly define a planar algebra 
which will be identified in section \ref{main} to be the one
corresponding to the {\em group-type subfactor} $P^H \subset P
\rtimes K$ of \cite{BH}.   

Let $G$ be a group generated by two of its finite subgroups $H$ and
$K$ and $e$ denote the identity element of $G$. Let us define
\[S_n = 
\left\{
\begin{array}{ll}
\{e\} & \text{if } n=0\\
\underbrace{K \times H \times K \times H \times \cdots}_{(n
 ~~ \text{factors})} & \text{if } n \geq 1   
\end{array}
\right.
\]
\[ S = \coprod_{n \geq 0} S_n \]
\[L_n = 
\left\{ 
\begin{array}{ll}
K, & \text{if $n$ is even} \\ H, & \text{otherwise}\\
\end{array} \right. \]
Let $\mu: S \rightarrow G$ be the multiplication map.
With the above notation, we are ready to define the planar algebra 
but first we would need some terminology.
\subsection*{Terminology:}
By a {\em face} in an unlabelled tangle $T$, we will mean a connected
component of
$D_0 \setminus [(\cup_{i=1}^b D_i) \cup \mathcal{S}]$ where $D_0$ is
the external disc, $D_i$ is the $i$-th internal disc for $i= 0, 1, 2
\cdots b$ and $\mathcal{S}$ is the set of strings of (an element in
the isotopy class of) $T$. By an {\em opening} in a tangle, 
we will mean the subset of the boundary of a disc lying  
between two consecutive marked points.
An opening will be called {\em internal} (resp., {\em external}) 
if it is a subset of the boundary of the internal (resp., external)
disc. Note that the boundary of a generic face may be disconnected
due to the presence of loops or networks inside it 
(see Figure \ref{face}). The set of connected components of the 
boundary of each face will have a single {\em outer component} 
and several (possibly none) {\em inner component(s)}.
\psfrag{d1}{$D_{i_1}$}
\psfrag{d2}{$D_{i_2}$}
\psfrag{d3}{$D_{i_3}$}
\psfrag{d4}{$D_{i_4}$}
\psfrag{d5}{$D_{i_5}$}
\psfrag{d6}{$D_{i_6}$}
\begin{figure}[h]
\begin{center}  
      \epsfig{file=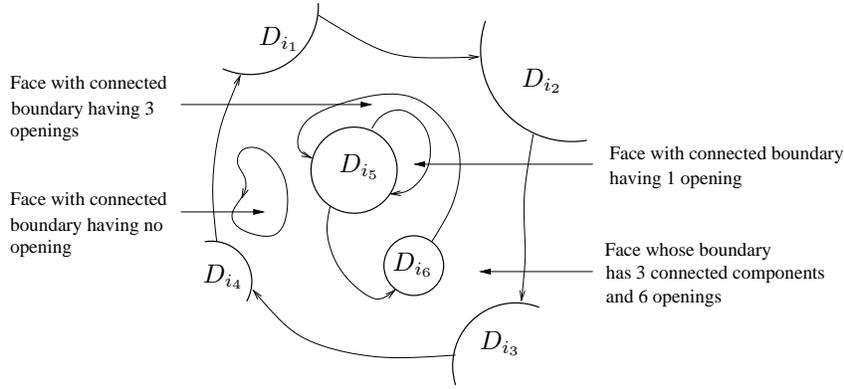, scale=0.5}
\caption{Example of faces in a tangle}
\label{face}
    \end{center}
\end{figure}
\begin{definition}\label{state}
A {\em state} $f$ on a tangle $T$ is a function 
$f: \{\text{all openings in } T\} \rightarrow H \coprod K$  such that
following holds:
\vskip 1em
\begin{itemize}
\item[(i)]
$f(\alpha) \in  \left \{ 
\begin{array}{ll} 
K, & \text{ if the face containing $\alpha$ is shaded,}\\ 
H, & \text{ otherwise.} \
\end{array}
\right. $
\vskip 1em
\item[(ii)]
{\em Triviality on the outer component of the boundary of a face}:

\noindent Let $\alpha_1, \alpha_2, \cdots
\alpha_m$ be the openings on the outer component
of the boundary of a face counted clockwise,
then we must have
\[ f(\alpha_1)^{\eta_1} f(\alpha_2)^{\eta_2} \cdots 
f(\alpha_m)^{\eta_m} = e\]
where $\eta_i = \left \{ 
\begin{array}{ll} 
+1, & \text{if $\alpha_i$ is an external opening,}\\
-1, & \text{otherwise.}
\end{array} 
\right.$
\vskip 1em
\item[(iii)]
{\em Triviality on internal discs}:

\noindent $f$ induces a map 
$\partial f : \{D_0, D_1, \cdots D_b\} \rightarrow \coprod_{n \geq
  0}S_n$
defined by 
\[ \partial f (D_i) = (f(\alpha^{(i)}_1), f(\alpha^{(i)}_2), 
\cdots f(\alpha^{(i)}_{2n_i})) \in S_{2n_i} \]
where $\alpha^{(i)}_1, \alpha^{(i)}_2, \cdots \alpha^{(i)}_{2n_i}$ 
are consecutive openings counted clockwise such that $\alpha^{(i)}_1$ 
is the opening between the first and the second marked points of 
$\partial D_i$.

\noindent We demand that $\mu(\partial f(D_i)) =\; e$, for all $i = 1,
2, \cdots b$.
\end{itemize}
\end{definition}

Note that triviality on the external disc and every boundary
component of every face follows (see Remark \ref{consequence} 
for a proof of this fact).

The above definition also applies to networks (positively or
negatively oriented).
\psfrag{a1}{$a_1$}
\psfrag{a2}{$a_2$}
\psfrag{a3}{$a_3$}
\psfrag{a4}{$a_4$}
\psfrag{a5}{$a_5$}
\psfrag{a6}{$a_6$}
\psfrag{b1}{$b_1$}
\psfrag{b2}{$b_2$}
\psfrag{b3}{$b_3$}
\psfrag{b4}{$b_4$}
\psfrag{c1}{$c_1$}
\psfrag{c2}{$c_2$}
\psfrag{c3}{$c_3$}
\psfrag{c4}{$c_4$}
\psfrag{d1}{$d_1$}
\psfrag{d2}{$d_2$}
\psfrag{d3}{$d_3$}
\psfrag{d4}{$d_4$}
\psfrag{d5}{$d_5$}
\psfrag{d6}{$d_6$}
\psfrag{d7}{$d_7$}
\psfrag{d8}{$d_8$}
\psfrag{f1}{$f_1$}
\psfrag{f2}{$f_2$}
\begin{figure}
\begin{center}  
      \epsfig{file=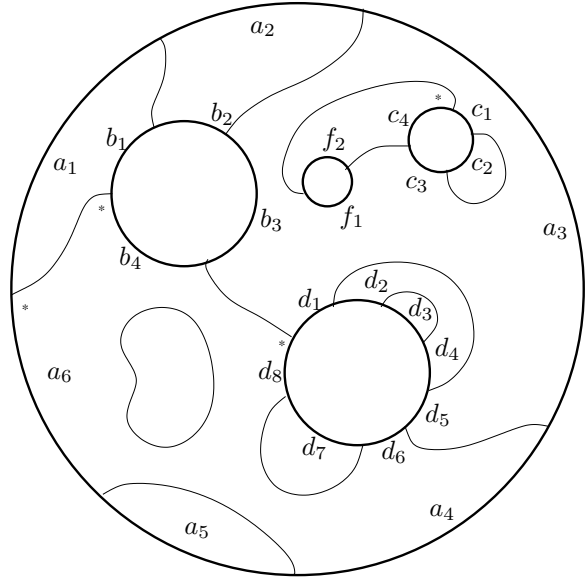, scale=0.5}
\caption{An example of a state on an unlabelled tangle}
\label{statepic}
    \end{center}
\end{figure}

Figure \ref{statepic} illustrates conditions (ii) and  
(iii) of Definition \ref{state}; triviality on internal 
discs give  
$b_1 b_2 b_3 b_4  =  
c_1 c_2 c_3 c_4  =   
d_1 d_2 d_3 d_4 d_5 d_6 d_7 d_8 =  
f_1 f_2 = e$
and triviality on the outer component of the boundaries of faces give 
$a_1 b_1^{-1} =  
a_2 b_2^{-1} =  
a_3 d_5^{-1}d_1^{-1} b_3^{-1} =
c_4^{-1} f_2^{-1} =  
c_2 =  
d_2^{-1} d_4^{-1} = 
d_3 =  
a_4 a_6 b_4^{-1} d_8^{-1} d_6^{-1} = 
d_7 =  
a_5 = e$. 
Note that the above relations also imply  
$a_1 a_2 a_3 a_4 a_5 a_6 = e$, and $c_1^{-1} c_3^{-1} f_1^{-1} = e$.

As the computation involving Figure \ref{statepic} suggests, 
triviality on inner boundaries of a face and on the external disc 
are actually consequences of the definition of a state. 
This is made precise in the following remark.
\begin{remark}\label{consequence}
Let $f$ be a state on a tangle or a network. Then the
following conditions hold: 
\begin{itemize}
\item[$\text{(ii)}^\prime$] {\em Triviality on every inner
  component of the boundary of a face}:

\noindent For every inner component of the boundary of a face with
openings $\alpha_1$, $\alpha_2$, $\dots$, $\alpha_m$ counted clockwise, we
have $f(\alpha_1) f(\alpha_2) \cdots f(\alpha_m) = e$.
\item[$\text{(iii)}^\prime$] {\em Triviality on the external disc}:

\noindent This is just a restatement of condition (iii) in Definition
\ref{state} applied to the external disc only if its color is greater
than zero.
\end{itemize}
\end{remark}

We prove the above remark using planar graphs. Without loss of
generality, we may start with a tangle or a network which is
{\em connected}. By a network or a tangle (with non-zero color on its
external disc) being connected, we mean that the union of the
boundaries of all the
discs and strings is a connected set; a $0$-tangle is said to be
connected if the network obtained after removing its external disc is
connected. To each tangle or network, we associate a planar graph
whose vertex set is the set of all marked points on the internal and
external discs, and the edges are the openings and the strings. Further,
we make this graph a directed one in such a way that the directions on
the edges arising from the openings are induced by clockwise
orientation on the boundary of the discs, and the remaining edges
(coming from the strings) are free to have any direction. Any state
$f$ assigns group elements to edges arising from openings; we label
each of the remaining
edges by $e$ and that is why we did not put
any restriction on the direction of such edges. Note that the
definition of the state implies the following condition on the group
labelled planar directed graph:

\noindent {\em Triviality on the boundary each bounded
  face of the graph}\footnote{Since we started with a connected tangle
  or network, the associated planar
graph will be connected; in particular, boundary of each face of the
graph will be connected}: 

\noindent If $g_1, g_2, \dots, g_m$ are the group elements assigned
to consecutive edges around a face read clockwise, then
\[g^{\eta_1}_1 g^{\eta_2}_2 \cdots  g^{\eta_m}_m = e\]
where $\eta_i = \left \{ 
\begin{array}{ll} 
+1, & \text{if $i$-th edge induces clockwise orientation in the
  face,}\\ -1, & \text{otherwise.}
\end{array} 
\right.$

\noindent To establish Remark \ref{consequence}, it is enough to
prove {\em triviality on the boundary of the unbounded face}. To see
this, we first consider a pair of bounded faces which have at least
one vertex or edge in common. If this pair is considered as a separate
graph, then using triviality on each face, it is easy to check
triviality on the boundary of the unbounded face of this pair. One can
then use this fact inductively to deduce the desired result for the
whole graph.
\vskip 1em
Next, we give the definition of the planar
algebra. Let the set of states on a tangle $T$ be denoted by
${\mathcal S}(T)$.
\subsection*{The vector spaces:} For $n \geq 0$, define $P_n = 
\C \{\ul{s} \in S_{2n}:\mu(\ul{s}) = e \}$.
\subsection*{Action of tangles:}
Let $T$ be an unlabelled tangle with (possibly zero) internal disc(s)
$D_1, D_2, \cdots , D_b$ and external disc $D_0$ where the color of
$D_i$ is $n_i$. Then $T$ defines a multilinear map, denoted by 
$Z_T : P_{n_1} \times P_{n_2} \cdots \times P_{n_b} \rightarrow
P_{n_0}$. We will define $Z_T(\ul{s_1}, \ul{s_2}, 
\cdots \ul{s_b}) \in P_{n_0}$, where $\ul{s_i} \in
S_{2n_i}$ such that $\mu(\ul{s_i}) = e$. 
Infact, we will just prescribe the coefficient of $\ul{s_0} \in
S_{2n_0}$ (such that $\mu(\ul{s_0}) = e$) in the expansion of 
$Z_T(\ul{s_1}, \ul{s_2}, \cdots \ul{s_b})$ in terms of the
canonical basis mentioned above.

We choose and fix a representative in the isotopy class of 
$T$ and call it the {\em standard form} of $T$, denoted by 
$\tilde{T}$. It is assumed to satisfy the following properties: 
\begin{itemize}
\item
$\tilde{T}$ is in rectangular form - meaning that - all of its 
discs are replaced by boxes and it is placed in $\mathbb{R}^2$ in such a
way that the boundaries of the boxes are parallel to the co-ordinate axes. 
\item
The first marked point on the boundary of each box is on the top left
corner.
\item
The collection of strings have finitely many local maxima or minima. 
\item
The external box can be sliced by horizontal lines in such a way that 
each maxima, minima, internal box is in a different slice. 
\end{itemize}
To each local maximum or minimum of a string with end-points, we
assign
a scalar according to Figure  \ref{assign}.
\psfrag{4throot}{$\sqrt[4]{\frac{\abs{H}}{\abs{K}}}$} 
\psfrag{inv4throot}{$\sqrt[4]{\frac{\abs{K}}{\abs{H}}}$} 
\begin{figure}
\begin{center}  
      \epsfig{file=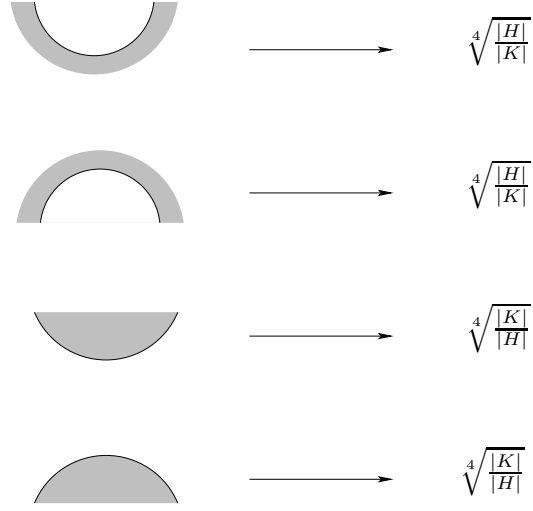, scale=0.5}
\caption{Assigning scalars to maxima or minima}       
\label{assign}
    \end{center}
\end{figure}
Let $p(T)$ denote the product of all scalars arising from the local
maxima or minima and $n_+ (T)$ (resp., $n_- (T)$) be the number of
non-empty connected
positively (resp., negatively) oriented network(s) in the tangle
$\tilde{T}$. Then, the coefficient $\langle Z_T(\ul{s_1}, \ul{s_2}, 
\cdots \ul{s_b}) | \ul{s_0} \rangle$ of $\ul{s_0}$ in $Z_T(\ul{s_1},
\ul{s_2}, \cdots \ul{s_b})$ is given by:
\[ p(T) \; \abs{H}^{n_+ (T)} \abs{K}^{n_- (T)} \; \abs{\{f \in
  {\mathcal S}(T) : \partial f(D_i) = \ul{s_i} \text{ for all } i =
  0,1, \cdots b\}} \]
Note that there could be several standard form representatives for
  $T$.
However, one standard form representative for $T$ can be transformed
into another by a finite sequence of moves of the following three
types:
\vskip 1em
\begin{itemize}
\item[I.] Horizontal or vertical sliding of boxes,
\item[II.] Wiggling of the strings,
\item[III.] Rotation of an internal box by a multiple of
  $360^{\circ}$.
\end{itemize}
\vskip 1em
It is easy to check that the above three moves do not alter the 
number of connected networks and keeps 
$\abs{ \{f \in {\mathcal S}(T) : \partial f(D_i) = \ul{s_i} \text{ for
all } i = 0,1, \cdots b\}}$ unaltered. So, it remains to show that
$p(T)$ is unchanged under the moves. Type I moves are the easiest because
they do not generate any new local maxima or minima. In each of the
moves of type II and III, there arises pair(s) of local maximum and
minimum in such a way that the two scalars assigned to each pair are
inverses of each other; as a result, $p(T)$ remains unchanged.
\subsection*{Action of the tangles preserve composition:} For $S$ an
$n_0$-tangle with internal discs $D_1 , D_2, \cdots , D_b$ having 
colors $n_1, n_2, \cdots, n_b$ respectively, and $T$ an $n_j$-tangle
for some $j \in \{1,2, \cdots b\}$, we would like to show
$Z_{S \circ_{D_j} T} = Z_S \circ (id_{P_{n_1}} \times \cdots 
Z_T  \cdots \times id_{P_{n_b}} )$.
For this, we first identify a set $\mathcal{E}$ of annular tangles
(with the distinguished internal disc as $D_1$) which we call {\em
  elementary tangles}, namely,
\psfrag{1}{$1$} \psfrag{2}{$2$}
\psfrag{n}{$n$} \psfrag{n-1}{${n-1}$} \psfrag{n+1}{$n+1$}
\psfrag{i}{$i$} \psfrag{i+1}{${i\!+\!1}$} 
\psfrag{D0}{$D_0$} \psfrag{D1}{$D_1$} \psfrag{D2}{$D_2$}
\psfrag{D0'}{$D^{\prime}_0$} 
\psfrag{a1}{$a_1$} \psfrag{ak}{$a_k$}
\psfrag{b1}{$b_1$} \psfrag{bl}{$b_l$}
\psfrag{E}{$E$} \psfrag{T}{$T$}
\psfrag{p strings}{$p$ strings} 
\psfrag{q strings}{$q$ strings} 
\psfrag{r strings}{$r$ strings} 
\psfrag{(i)}{(i) Capping tangles:}
\psfrag{(ii)}{(ii) Cap inclusion tangles:}
\psfrag{(iii)}{(iii) Left inclusion tangles:}
\psfrag{(iv)}{(iv) Disc inclusion tangle:}
\psfrag{(v)}{$\text{(iv)}^\prime$ Disc inclusion tangle:}
\psfrag{c}{$\cdots$}
\psfrag{*}{$*$}
\psfrag{kh}{$(k,h)$}
\begin{flushleft}
\epsfig{file=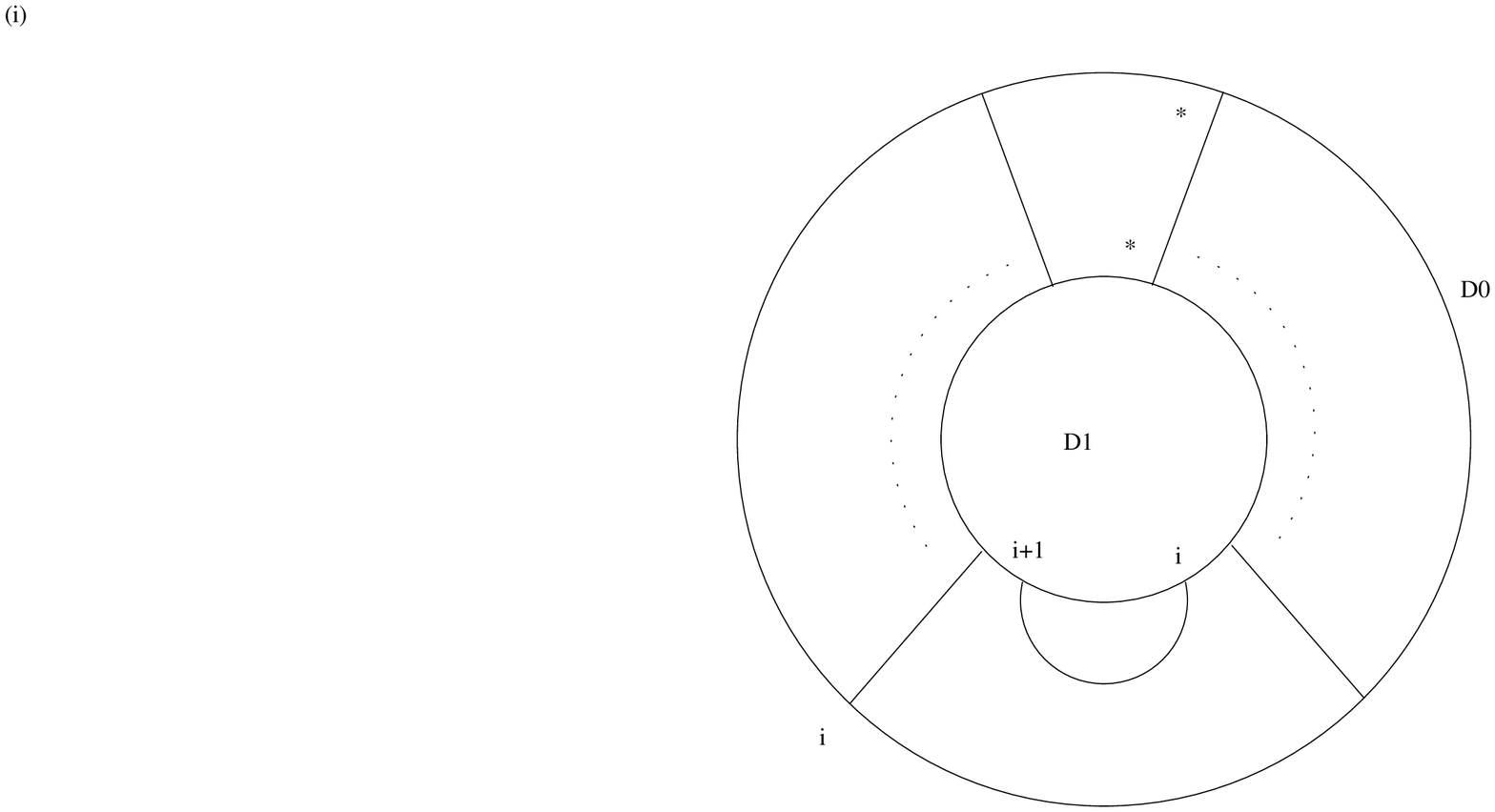, scale=0.35}
\end{flushleft}
\noindent with $col(D_1) = n \geq 1,\; col(D_0) = n-1$, and 
$1 \leq i \leq 2n-1$ ($col(D_1) = 1$ just means that there are 
no strings connecting the internal disc to the external disc).
\begin{flushleft}  
      \epsfig{file=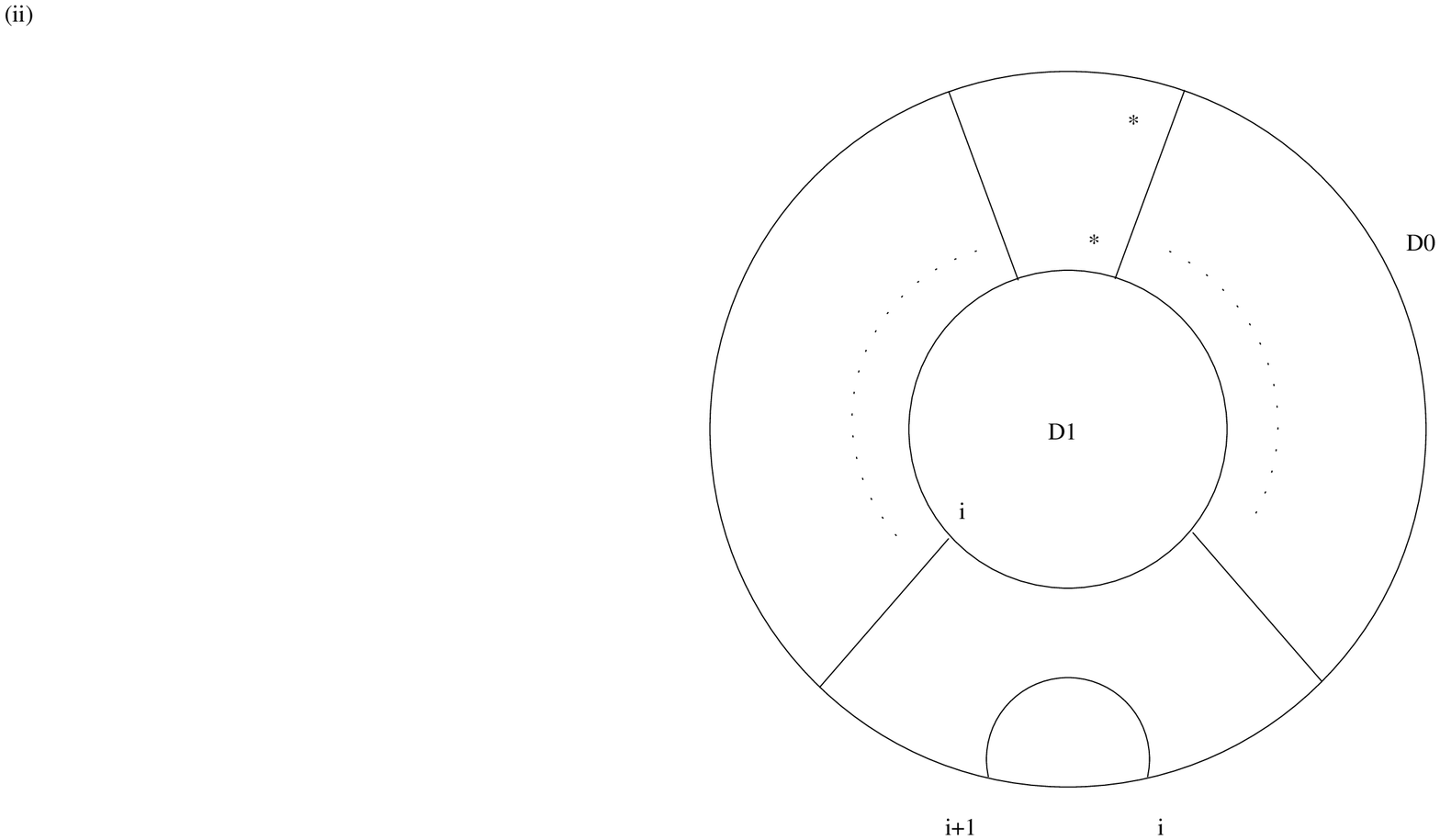, scale=0.35}
\end{flushleft}
with $col(D_1) = n \geq 0,\; col(D_0) = n+1$ , and $1 \leq i \leq
2n+1$.\\

If $i=1$ in (i) or (ii) then the pictures should be interpreted 
as simply not having the bunch of $(i-1)$ straight strings joining 
the marked points of the internal disc and the corresponding points 
of the external disc.
\begin{flushleft}  
      \epsfig{file=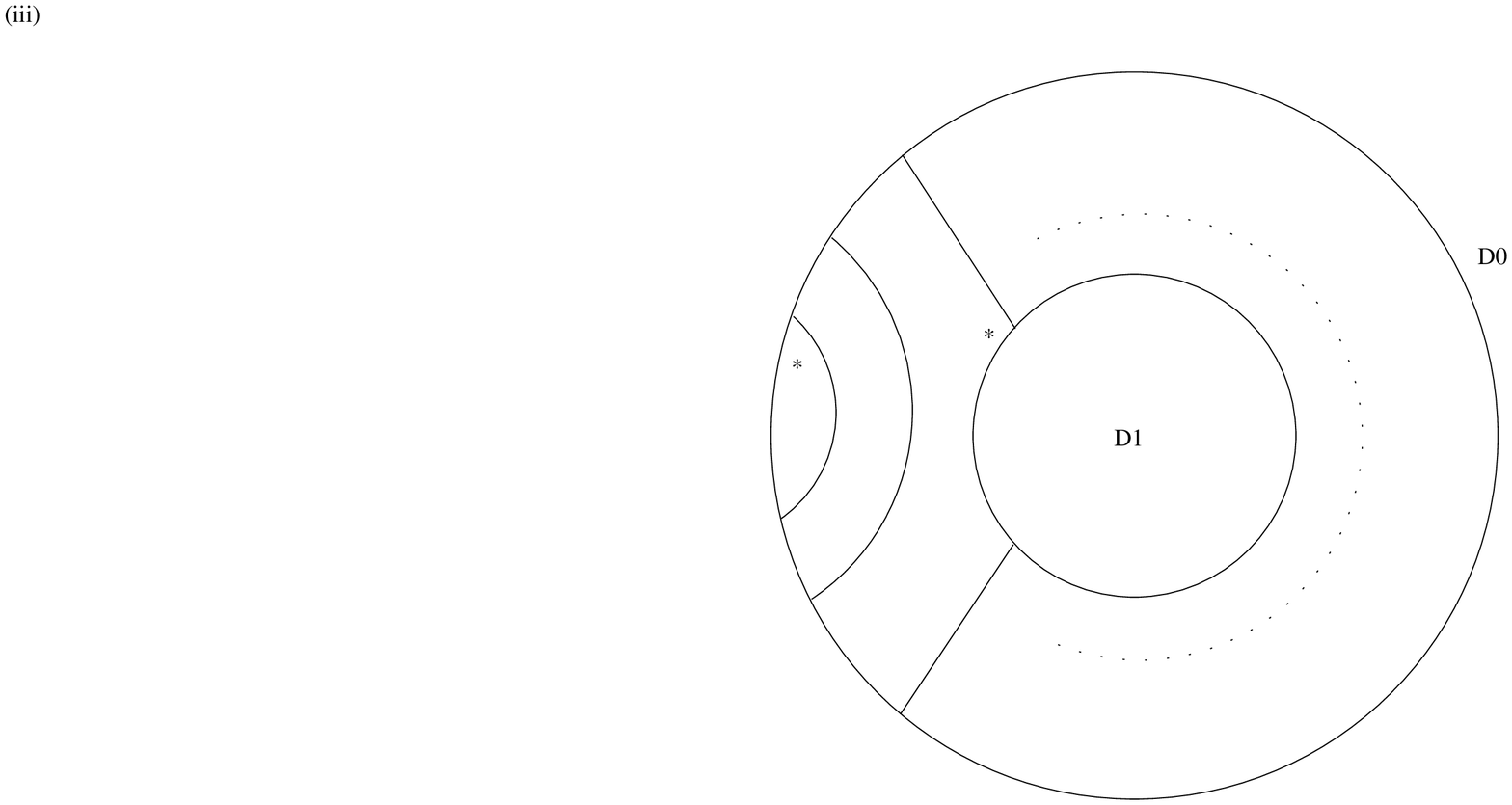, scale=0.35}
\end{flushleft}
with $col(D_1) = n \geq 0,\; col(D_0) = n+2$.
\psfrag{DIT}{ }
\begin{flushleft}  
      \epsfig{file=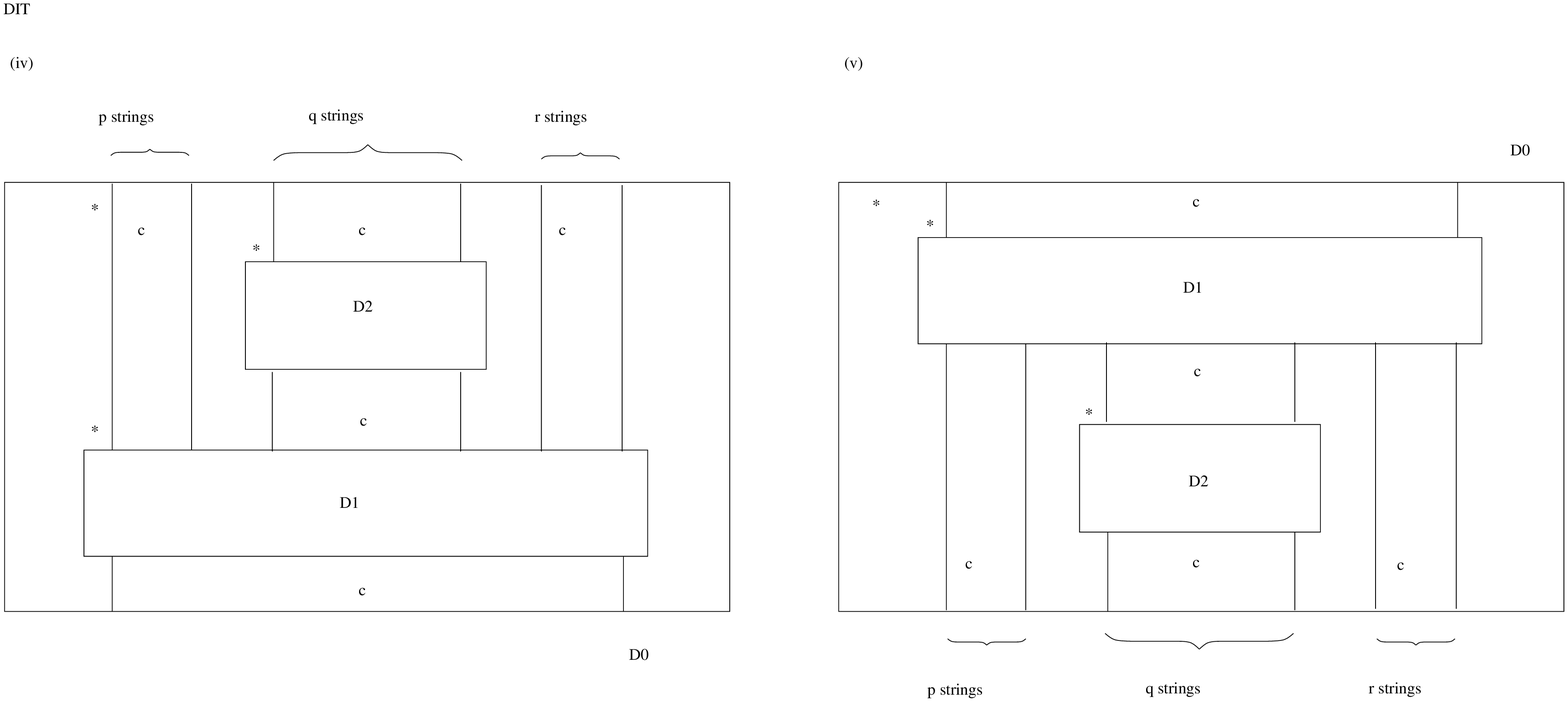, scale=0.34}
\end{flushleft}
\noindent with $col(D_1) = col(D_0) = n \geq 0,\; col(D_2) = q,\;
\text{where}\;
p \geq 0, q \geq 0, r \geq 0$~ such that ~$p+q+r = n$ and $p$ is even.
\vskip 1em
Note that any annular tangle can be expressed as a
composition of the elementary tangles. To see this, express the
annular tangle in standard form and cut it into horizontal strips
each of which contains at most one internal disc or one local maxima 
or minima. Now, the strip containing the distinguished
internal disc inside the annular tangle can be obtained by 
composition of elementary tangles of type (iii) and type (ii) 
(more specifically, the inclusion tangles); one can then glue the 
other strips consecutively one after the other along the lines of 
cutting to get back the original tangle. Each such gluing operation 
is given by composition of an elementary tangle of types (i), 
(ii), (iv) or
$\text{(iv)}^\prime$. So, to prove that the action of the tangles
preserve composition, it is enough to prove $Z_{E \circ_{D_1} T} = Z_E
\circ Z_T$ (resp., 
$Z_{E \circ_{D_1} T} = Z_E \circ (Z_T \times id_{P_q} )$)
for any tangle $T$ and any $E \in \mathcal{E}$ of type (i),
(ii) or (iii) (resp. (iv) or $\text{(iv)}^\prime$) whenever the
composition makes sense.
\vskip 1em
We fix an $n$-tangle $T$ with internal discs $D^{\prime}_1 ,
D^{\prime}_2, \cdots , D^{\prime}_b$ with
colors $n_1, n_2, \cdots, n_b$ respectively, and an $n_0$-tangle $E
\in \mathcal{E}$ such that both $T$ and $E$ are in standard forms and
color of 
$D_1$ in $E$ is $n$. Let us consider the standard form on $E
\circ_{D_1} T$ induced by the standard forms of $E$ and $T$. Our goal
is to show:
\[ \langle Z_{E \circ_{D_1} T} (\ul{s_1},
\ul{s_2}, \cdots \ul{s_b}) | \ul{s_0} \rangle \; =
\sum_{\substack {\ul{s} \; \in \; S_{2n} \\ \text{s.t. }
    \mu(\ul{s})=e}}
\langle Z_E (\ul{s}) | \ul{s_0} \rangle \;
\langle Z_T (\ul{s_1},
\ul{s_2}, \cdots \ul{s_b}) | \ul{s} \rangle\]
if $E$ is of type (i), (ii) or (iii), and
\[\langle Z_{E \circ_{D_1} T} (\ul{s_1},
\ul{s_2}, \cdots \ul{s_b} , \ul{t}) | \ul{s_0} \rangle \; =
\sum_{\substack {\ul{s} \; \in \; S_{2n} \\ \text{s.t. }
    \mu(\ul{s})=e}}
\langle Z_E (\ul{s} , \ul{t}) | \ul{s_0} \rangle \;
\langle Z_T (\ul{s_1},
\ul{s_2}, \cdots \ul{s_b}) | \ul{s} \rangle \]
if $E$ is of type (iv) or $\text{(iv)}^\prime$ where $\ul{s_j} \in
S_{2n_j}$
for $0 \leq j \leq b$ and $\ul{t} \in S_q$ such that $\mu
(\ul{s_j})=e=
\mu (\ul{t})$ for all $j$. 
An interesting situation arises when we pick elementary tangles of
type
(i), since composition of tangles in this case may lead to a change in
the number of connected networks. The reasoning in the other cases is 
either similar or straightforward.   
\vskip 1em
For $E$ being type (i) elementary tangle, the above equation is
equivalent to:
\[ \begin{array}{ll} 
& p(E \circ_{D_1} T) \, \abs{H}^{n_+ (E \circ_{D_1} T)}
\, \abs{K}^{n_- (E \circ_{D_1} T)}
\left |
\left \{
f \in {\mathcal S}(E \circ_{D_1} T)
\left |
\begin{array}{c}
\partial f(D^{\prime}_j) = \ul{s_j}\\
\text{ for } 1 \leq j \leq b,\\
\partial f(D_0) = \ul{s_0}
\end{array}
\right .
\right \} \right | \\
&\\
= & p(E) \; p(T) \; \abs{H}^{n_+ (E) + n_+ (T)}
\; \abs{K}^{n_- (E) + n_- (T)}\\
& \cdot \sum_{\substack {\ul{s} \; \in \; S_{2n} \\ \text{s.t. }
\mu(\ul{s})=e}}
\left |
\left \{
f \in {\mathcal S}(E)
\left |
\begin{array}{l}
\partial f(D_1) = \ul{s}\\
\partial f(D_0) = \ul{s_0}
\end{array}
\right .
\right \}
\right | \,
\left |
\left \{
f \in {\mathcal S}(T)
\left |
\begin{array}{c}
\partial f(D^{\prime}_j) = \ul{s_j}\\
\text{ for } 1 \leq j \leq b,\\
\partial f(D^{\prime}_0) = \ul{s}
\end{array}
\right .
\right \}
\right |
\end{array} \]
where $D^{\prime}_0$ denotes the external disc of $T$. First, observe
that
$p(E \circ_{D_1} T) = p(E) \, p(T)$. To show the equality of the 
remaining scalars, we consider the following two cases.
\vskip 1em
\noindent {\bf Case 1}: The string which connects the $i$-th and the
$(i+1)$-th points on $D_1$ in $E$, does not produce any new network in
$E \circ_{D_1} T$ other than those that are already present in $T$.
\vskip 1em
Clearly, $n_{\epsilon} (E \circ_{D_1} T) = n_{\epsilon}(T)$
(since no new network appears in
$E \circ_{D_1} T$) and $n_{\epsilon} (E) = 0$ for
$\epsilon \in \{+, -\}$.

A typical example of
such a case can be viewed in the following picture where we label the
openings on the internal discs of $T$ by group elements coming from
the coordinates of $\ul{s_j}$ for $1 \leq j \leq b$, and the openings
on $D_0$ of $E$ by coordinates of $\ul{s_0} = (g_1,g_2, \cdots,
g_{2n-2})$. 
\psfrag{gi}{$g_i$}
\psfrag{gi-1}{$g_{i-1}$}
\psfrag{gi-2}{$g_{i-2}$}
\begin{figure}[h]
      \epsfig{file=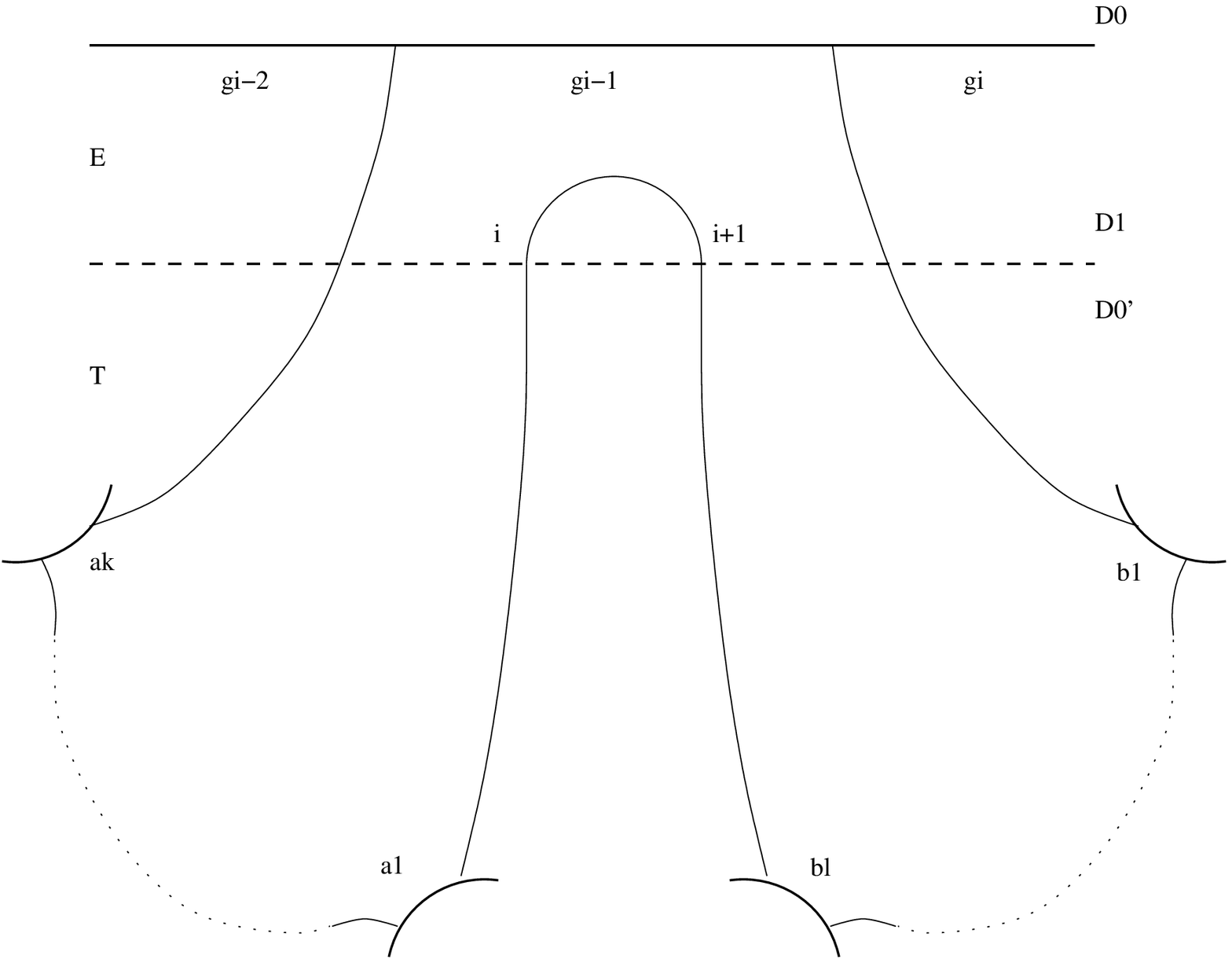, scale=0.4}
\caption{}
\label{composition1}
\end{figure}\\
\noindent Using triviality on the boundary of faces in the definition
of a state, we get
\[
\begin{array}{l}
\left | \left \{
f \in {\mathcal S}(E)
\left |
\begin{array}{c}
\partial f(D_1) = \ul{s}\\
\partial f(D_0) = \ul{s_0}
\end{array}
\right .
\right \} \right | \neq 0 \Longleftrightarrow
\left | \left \{
f \in {\mathcal S}(E)
\left |
\begin{array}{c}
\partial f(D_1) = \ul{s}\\
\partial f(D_0) = \ul{s_0}
\end{array}
\right .
\right \} \right | = 1\\
\\
\Longleftrightarrow \ul{s} = (g_1, \cdots,
g_{i-2}, (g_{i-1} g), e , g^{-1}, g_i, \cdots, g_{2n-2}) \in S_{2n} \text{
  for some } g \in L_i.
\end{array}
\]
Define $\ul{s}^g = (g_1, \cdots,
g_{i-2}, (g_{i-1} g), e , g^{-1}, g_i, \cdots, g_{2n-2})$ for $g \in
L_i$. So, it is enough to check
\[
\left | \left \{
f \in {\mathcal S}(E \circ_{D_1} T)
\left |
\begin{array}{c}
\partial f(D^{\prime}_j) = \ul{s_j}\\
\text{ for } 1 \leq j \leq b,\\
\partial f(D_0) = \ul{s_0}
\end{array}
\right .
\right \} \right |
=
\sum_{g \in L_i}
\left | \left \{
f \in {\mathcal S}(T)
\left |
\begin{array}{c}
\partial f(D^{\prime}_j) = \ul{s_j}\\
\text{ for } 1 \leq j \leq b,\\
\partial f(D^{\prime}_0) = \ul{s}^g
\end{array}
\right .
\right \} \right |
\]
Carefully observing  Figure \ref{composition1} and using triviality on the
boundary of faces once again, we get
\[
\begin{array}{cl}
& \left | \left \{
f \in {\mathcal S}(E \circ_{D_1} T)
\left |
\begin{array}{c}
\partial f(D^{\prime}_j) = \ul{s_j}\\
\text{ for } 1 \leq j \leq b,\\
\partial f(D_0) = \ul{s_0}
\end{array}
\right .
\right \} \right | 
\neq 0 \text{, equivalently, equals to } 1\\
&\\
\Longrightarrow &
\left | \left \{
f \in {\mathcal S}(T)
\left |
\begin{array}{c}
\partial f(D^{\prime}_j) = \ul{s_j}\\
\text{ for } 1 \leq j \leq b,\\
\partial f(D^{\prime}_0) = \ul{s}^g
\end{array}
\right .
\right \} \right |
= \delta_{g \, , \, b^{\eta_1}_{1} b^{\eta_2}_{2} \cdots
  b^{\eta_l}_{l}}
\end{array}
\]
where $\eta_j = \pm 1$ according as the corresponding opening is
external or internal. Conversely, if 
$
\sum_{g \in L_i}
\left | \left \{
f \in {\mathcal S}(T)
\left |
\begin{array}{c}
\partial f(D^{\prime}_j) = \ul{s_j}\\
\text{ for } 1 \leq j \leq b,\\
\partial f(D^{\prime}_0) = \ul{s}^g
\end{array}
\right .
\right \} \right |
$ is non-zero,
then it has to equal to $1$ since $g$ must be $b^{\eta_1}_{1} b^{\eta_2}_{2}
\cdots b^{\eta_l}_{l}$ by triviality on the boundary of faces in $T$;
from the unique state on $T$ which makes the above sum non-zero,
one can easily induce a well-defined state on $E \circ_{D_1} T$, and hence
$\left | \left \{
f \in {\mathcal S}(E \circ_{D_1} T)
\left |
\begin{array}{c}
\partial f(D^{\prime}_j) = \ul{s_j}\\
\text{ for } 1 \leq j \leq b,\\
\partial f(D_0) = \ul{s_0}
\end{array}
\right .
\right \} \right | = 1$.
This completes the proof of Case 1.
\vskip 1em
\noindent {\bf Case 2}: The string which connects the $i$-th and the
$(i+1)$-th points on $D_1$ in $E$, produces a new (connected) network in
$E \circ_{D_1} T$ other than those that are already present in $T$.
\vskip 1em
First, let us assume that the new network is
positively oriented, equivalently, $i$ is odd.
Clearly, $n_- (E \circ_{D_1} T) = n_-(T)$ (since no
new negatively oriented network appears in $E \circ_{D_1} T$), and
$n_+ (E \circ_{D_1} T) = n_+(T) + 1$. 

Further, assume that $col(D_0) \geq 1$.
In this case, $n_\epsilon(E)=0$ for $\epsilon \in \{+,-\}$. 
A typical example of this case 
can be viewed in the following picture where we label the
openings on the internal discs of $T$ by group elements coming from
the coordinates of $\ul{s_j}$ for $1 \leq j \leq b$, and the openings
on $D_0$ of $E$ by coordinates of $\ul{s_0} = (g_1,g_2, \cdots,
g_{2n-2})$.
\begin{figure}[h]
\begin{center}  
      \epsfig{file=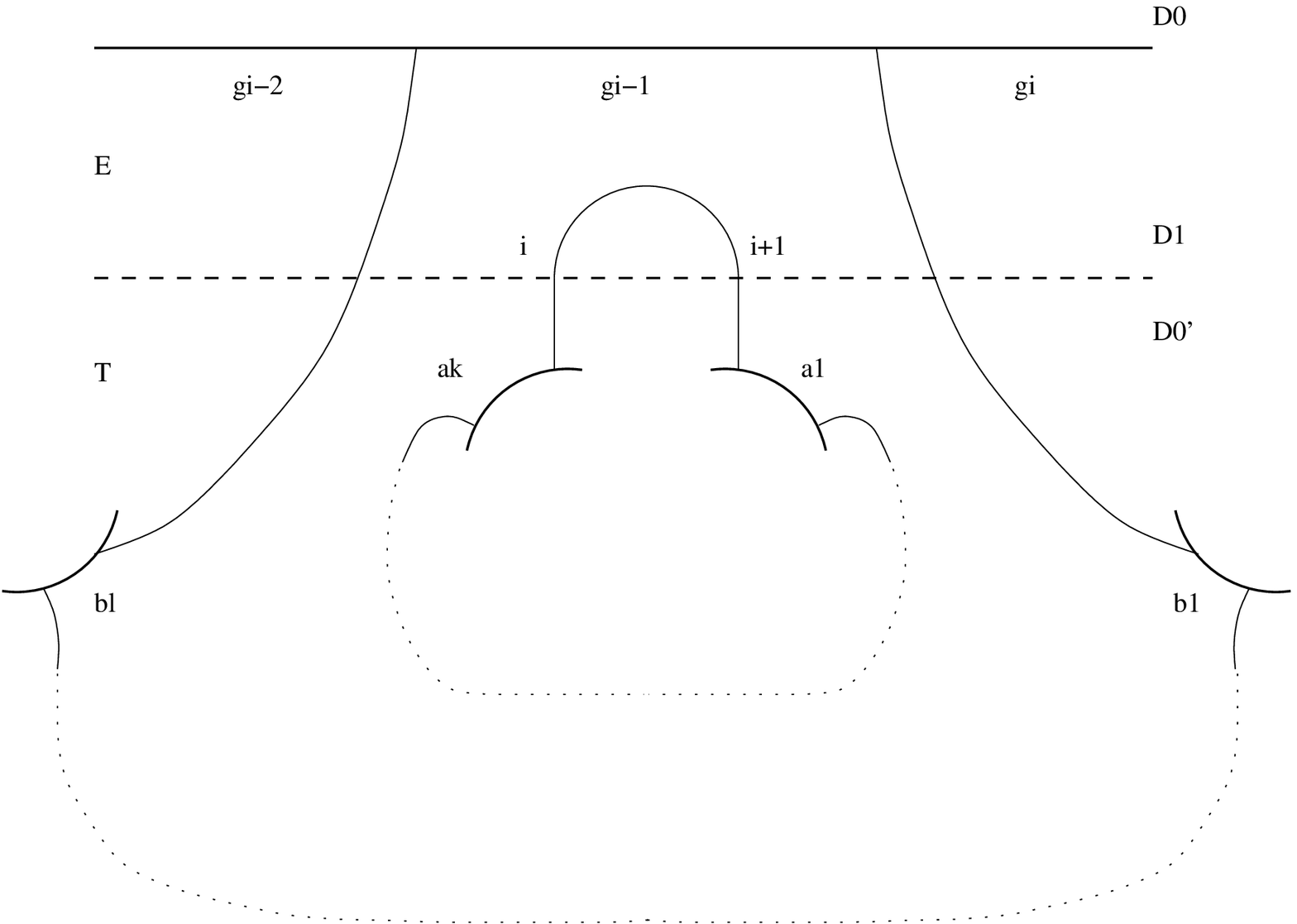, scale=0.4}
\caption{}
\label{composition2}
    \end{center}
\end{figure}
\noindent So, in this case, it is enough to check
\[
\abs{H} \left | \left \{
f \in {\mathcal S}(E \circ_{D_1} T)
\left |
\begin{array}{c}
\partial f(D^{\prime}_j) = \ul{s_j}\\
\text{ for } 1 \leq j \leq b,\\
\partial f(D_0) = \ul{s_0}
\end{array}
\right .
\right \} \right |
=
\sum_{g \in H}
\left | \left \{
f \in {\mathcal S}(T)
\left |
\begin{array}{c}
\partial f(D^{\prime}_j) = \ul{s_j}\\
\text{ for } 1 \leq j \leq b,\\
\partial f(D^{\prime}_0) = \ul{s}^g
\end{array}
\right .
\right \} \right |
\]
If 
$\left \{
f \in {\mathcal S}(E \circ_{D_1} T)
\left |
\begin{array}{c}
\partial f(D^{\prime}_j) = \ul{s_j}\\
\text{ for } 1 \leq j \leq b,\\
\partial f(D_0) = \ul{s_0}
\end{array}
\right .
\right \}$
is non-empty (equivalently, singleton), from Figure
\ref{composition2}, we have $a_1 a_2 \cdots a_k = e$ and
$g_{i-1} b^{\eta_1}_{1} b^{\eta_2}_{2} \cdots b^{\eta_l}_{l} = e$
where $\eta_j = \pm 1$ according as the corresponding
opening is external or internal. For any $g \in H$, define $f^g$ by
setting $\partial f(D^{\prime}_j) = \ul{s_j}$ for $1 \leq j \leq b$
and $\partial f(D^{\prime}_0) = \ul{s}^g$. To check whether $f^g$ is a
state, we consider the face in $T$ appearing in Figure
\ref{composition2};
triviality on the boundary of this face is given by the equation
$g^{-1} b^{\eta_1}_{1} b^{\eta_2}_{2} \cdots b^{\eta_l}_{l} (g_{i-1} g)
a_1 a_2 \cdots a_k = e$ which indeed holds. Triviality on all other
discs or faces are induced by the existence of the state on $E
\circ_{D_1} T$. Thus, $\left | \left \{
f \in {\mathcal S}(T)
\left |
\begin{array}{c}
\partial f(D^{\prime}_j) = \ul{s_j}\\
\text{ for } 1 \leq j \leq b,\\
\partial f(D^{\prime}_0) = \ul{s}^g
\end{array}
\right .
\right \} \right | = 1$ for all $g \in H$. Conversely, if the right
hand side is non-zero, then there exists $g \in H$ such that $f^g$ (defined
earlier) is a state on $T$. Analyzing Figure \ref{composition2}, we get
$g^{-1} b^{\eta_1}_{1} b^{\eta_2}_{2} \cdots b^{\eta_l}_{l} (g_{i-1} g)
a_k^{-1} \cdots a_2^{-1} a_1^{-1}= e$. 
Note that the opening between the $i$-th and
the $(i+1)$-th points of the disc $D^{\prime}_0$ of $T$, is assigned
$e$. Now, if we consider the network appearing in Figure \ref{composition2}
separately, then we have triviality on each of its internal faces and discs
(induced by $f^g$ being a state); by Remark \ref{consequence}
$\text{(ii)}^\prime$, we also have $a_1 a_2 \cdots a_k = e$
(triviality on the internal boundary of the external face of the
network). This implies 
$b^{\eta_1}_{1} b^{\eta_2}_{2} \cdots b^{\eta_l}_{l} g_{i-1} =
e$; as a result, $f^h$ is a state for every $h \in H$. So, if right hand
side is non-zero, then it has to be $\abs{H}$; moreover, we get
$b^{\eta_1}_{1} b^{\eta_2}_{2} \cdots b^{\eta_l}_{l} g_{i-1} = e$
which plays an important role in showing that 
$\left | \left \{
f \in {\mathcal S}(E \circ_{D_1} T)
\left |
\begin{array}{c}
\partial f(D^{\prime}_j) = \ul{s_j}\\
\text{ for } 1 \leq j \leq b,\\
\partial f(D_0) = \ul{s_0}
\end{array}
\right .
\right \}
\right | \neq 0$ (equivalently, equals to $1$).

This finishes the proof for the case where $i$ is odd.
For $i$ even, the proof is exactly similar, except that 
one has to interchange $\abs{H}$ and $\abs{K}$.
\vskip 1em
The subcase that deserves separate treatment is when 
$col(D_0) = 0$. In this case, $n_+(E)=1$ and 
$| \{ f \in {\mathcal S}(E) : \partial f(D_1) = \ul{s} \} | =
\delta_{\ul{s}, (e,e)}$. Therefore, it is enough to show
\[
\left | \left \{
f \in {\mathcal S}(E \circ_{D_1} T)
\left |
\begin{array}{c}
\partial f(D^{\prime}_j) = \ul{s_j}\\
\text{ for } 1 \leq j \leq b\\
\end{array}
\right .
\right \} \right |
= 
\left | \left \{
f \in {\mathcal S}(T)
\left |
\begin{array}{c}
\partial f(D^{\prime}_j) = \ul{s_j}\\
\text{ for } 1 \leq j \leq b,\\
\partial f(D^{\prime}_0) = (e,e)
\end{array}
\right .
\right \} \right |
\]

The proof of the equality of the two sides is similar 
and is left to the reader.

\vskip 1em
We now analyze the filtered $*$-algebra structure of $P$ and the
action of Jones projection tangles and conditional expectation tangles
which will be useful in section \ref{main} to show that $P$ is
isomorphic to the planar algebra arising from a group-type 
subfactor. We start with laying some notations. Define
\vskip 1em
\[\widetilde{S_n} =
\left\{
\begin{array}{ll}
\{e\} & \text{if } n=0\\
\underbrace{\cdots \times H \times K \times H \times K}_{(n factors)}
& \text{if } n \geq 1 
\end{array}
\right.
\]
\[T_n =
\left\{
\begin{array}{ll}
\{e\} & \text{if } n=0\\
\underbrace{H \times K \times H \times K \times \cdots}_{(n factors)}
& \text{if } n \geq 1 
\end{array}
\right.
\]
\vskip 1em
Define ${}^\sim : S_n \rightarrow \widetilde{S_n}$ by $(s_1 , s_2 , \cdots ,
  s_n)\widetilde{} = (s_n^{-1} , \cdots , s_2^{-1} , s_1^{-1})$ for
  $(s_1 , s_2 , \cdots , s_n) \in S_n$ and let ${}^- ~ :
  \widetilde{S_n} \rightarrow S_n$ denote its inverse.
\begin{remark}\label{structure}
We describe below the main structural features of the planar algebra $P$.
\begin{itemize}
\item[(i)] {\em Identity}: 
\[1_{P_n} = 
\left\{
\begin{array}{ll}
{e} & \text{if } n = 0\\
\sum_{\ul{s} \in S_{n-1}}
  (\ul{s} , e , \widetilde{\ul{s}}, e) & \text{if } n \geq 1
\end{array}
\right.
\]
\item[(ii)] {\em $*$-structure}: Define $*$ on $P$ by defining on the
  basis as
\[{\ul{s}}^{*} =
  ({s_{2n-1}}^{-1} , \cdots , {s_2}^{-1} , {s_1}^{-1} , {s_{2n}}^{-1})
  \]
where $\ul{s} =(s_1 , s_2 , \cdots , s_{2n}) \in S_{2n}$ such that
$\mu(\ul{s}) = e$ for $n \geq 1$, and then extend conjugate
linearly. Clearly, $*$ is an involution. One also needs to verify
whether the action of a tangle $T$ preserves $*$, that is, $Z_{T^*}
\circ (* \times \cdots \times *) = * \circ Z_T$; in particular,
$\langle Z_{T^*} (\ul{s^*_1}, \cdots, \ul{s^*_b}) | \ul{s^*_0} \rangle
= \overline{\langle Z_T (\ul{s_1}, \cdots, \ul{s_b}) | \ul{s_0}
  \rangle}$.  It is enough to check this equation for the cases when
$T$ has no internal disc or closed loops, and when $T$ is an
elementary tangle. The actual verification in each of these cases is
completely routine and is left to the reader.
\item[(iii)] {\em Multiplication}:
\[(\ul{a_1} , l_1 , \ul{b_1} , h_1) \cdot (\ul{a_2} , l_2 , \ul{b_2} ,
  h_2) = \delta_{\ul{b_1} ~ , ~ \widetilde{\ul{a_2}}} ~ (\ul{a_1} , l_1 l_2
  , \ul{b_2} , h_2 h_1)\]
where $\ul{a_i} \in S_{n-1}$, $\ul{b_i} \in \widetilde{S_{n-1}}$, $l_i
  \in L_{n-1}$, $h_i \in H$ such that $\mu(\ul{a_i} , l_i , \ul{b_i} ,
  h_i) = e$ for $i = 1 , 2$ and $n \geq 1$ (where we consider the
  elements $\ul{a_i}$ and $\ul{b_i}$ to be void in the case of $n = 1$),
\item[(iv)] {\em Inclusion}:
\[P_n \ni \ul{s} \mapsto \sum_{\substack{l_1 , ~ l_2 \in
    L_{n-1}\\\text{such that } l_1 l_2 = s_n}} (s_1 , s_2 , \cdots ,
    s_{n-1} , l_1 , e , l_2 , s_{n+1} , \cdots s_{2n}) \in P_{n+1}\]
where $\ul{s} =(s_1 , s_2 , \cdots , s_{2n}) \in S_{2n}$ such that
$\mu(\ul{s}) = e$ for $n \geq 1$,
\item[(v)] {\em Jones Projection Tangle}:\\ 
For $P_2$,  
\psfrag{=1st jonesproj}
{$ = \sqrt{\frac{\abs{K}}{\abs{H}}}
   \sum_{h \in H}(e, h, e, h^{-1})$}
\begin{flushleft}
      \epsfig{file=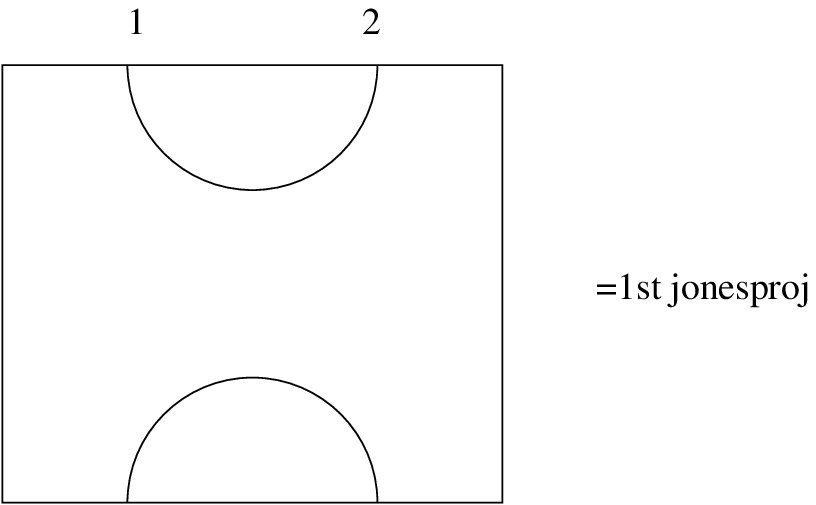, scale=0.5}
\end{flushleft}
and for $P_{n+1}$ for $n > 1$,
\psfrag{= formula for jonesproj}
{$= \sqrt{\frac{\abs{L_{n-1}}}{\abs{L_n}}}
     \sum_{\substack{\ul{s} \in S_{n-2}\\
          l_1, l_2, l_3 \in L_{n}\\
          \text{ s.t. } l_1 l_2 l_3 = e}} 
  (\ul{s}, l_1, e, l_2, e, l_3, \tilde{\ul{s}}, e)$}
\begin{flushleft}  
      \epsfig{file=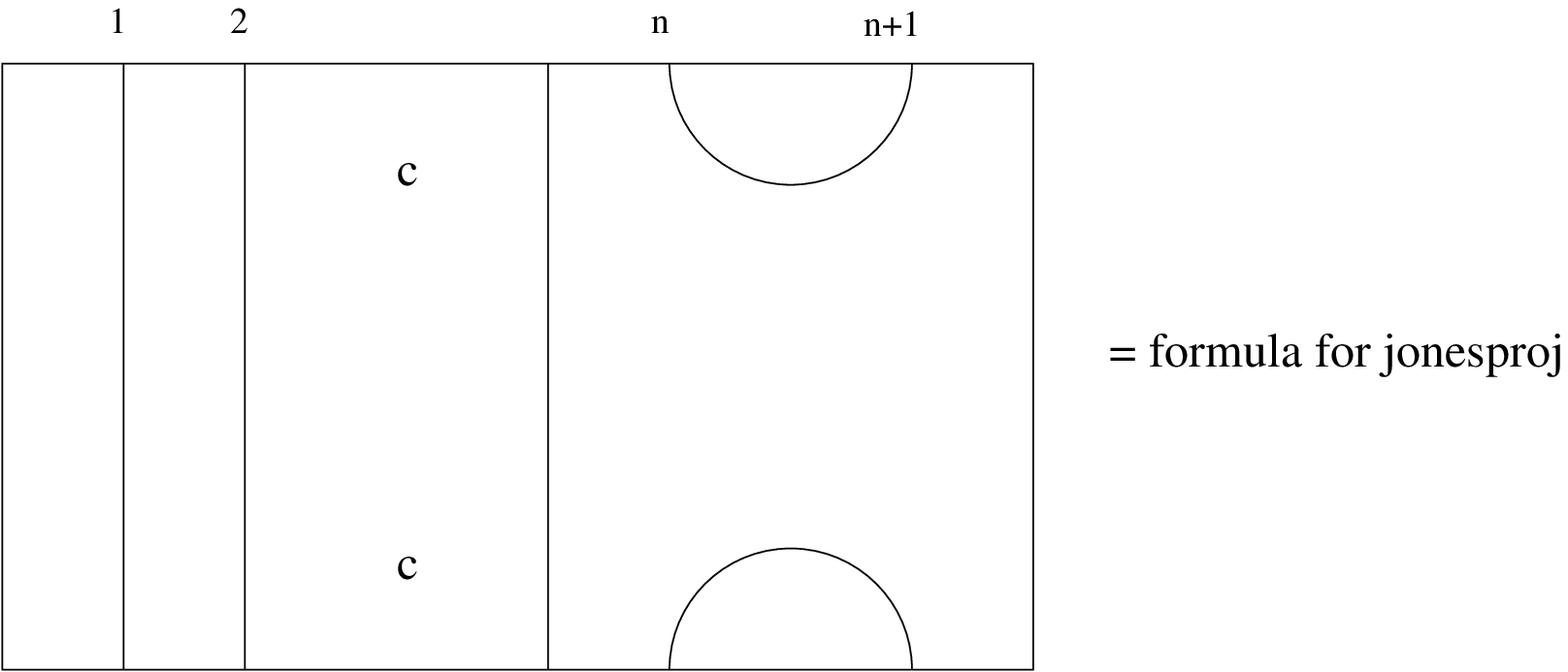, scale=0.4}
\end{flushleft}
\item[(vi)] {\em Conditional Expectation Tangle from $P_{n+1}$ onto $P_n$}:\\ 
For $n \geq 1$, let $\ul{s_1} \in S_{n-1}$, $\ul{s_2} \in \widetilde{S_{n-1}}$,
$m_1, m_2 \in L_{n-1}$, $l \in
L_n$, $h \in H$ such that $\mu(\ul{s_1},\, m_1 m_2,\, \ul{s_2},\, h) =
e$. Then, 
\psfrag{input}{$(\ul{s_1}, m_1, l, m_2, \ul{s_2}, h)$}
\psfrag{formula for condexp}
{$= \delta_{l,e}\;\sqrt{\frac{\abs{L_n}}{\abs{L_{n-1}}}}
\;(\ul{s_1},\, m_1 m_2,\, \ul{s_2},\, h)$}
\begin{flushleft}  
      \epsfig{file=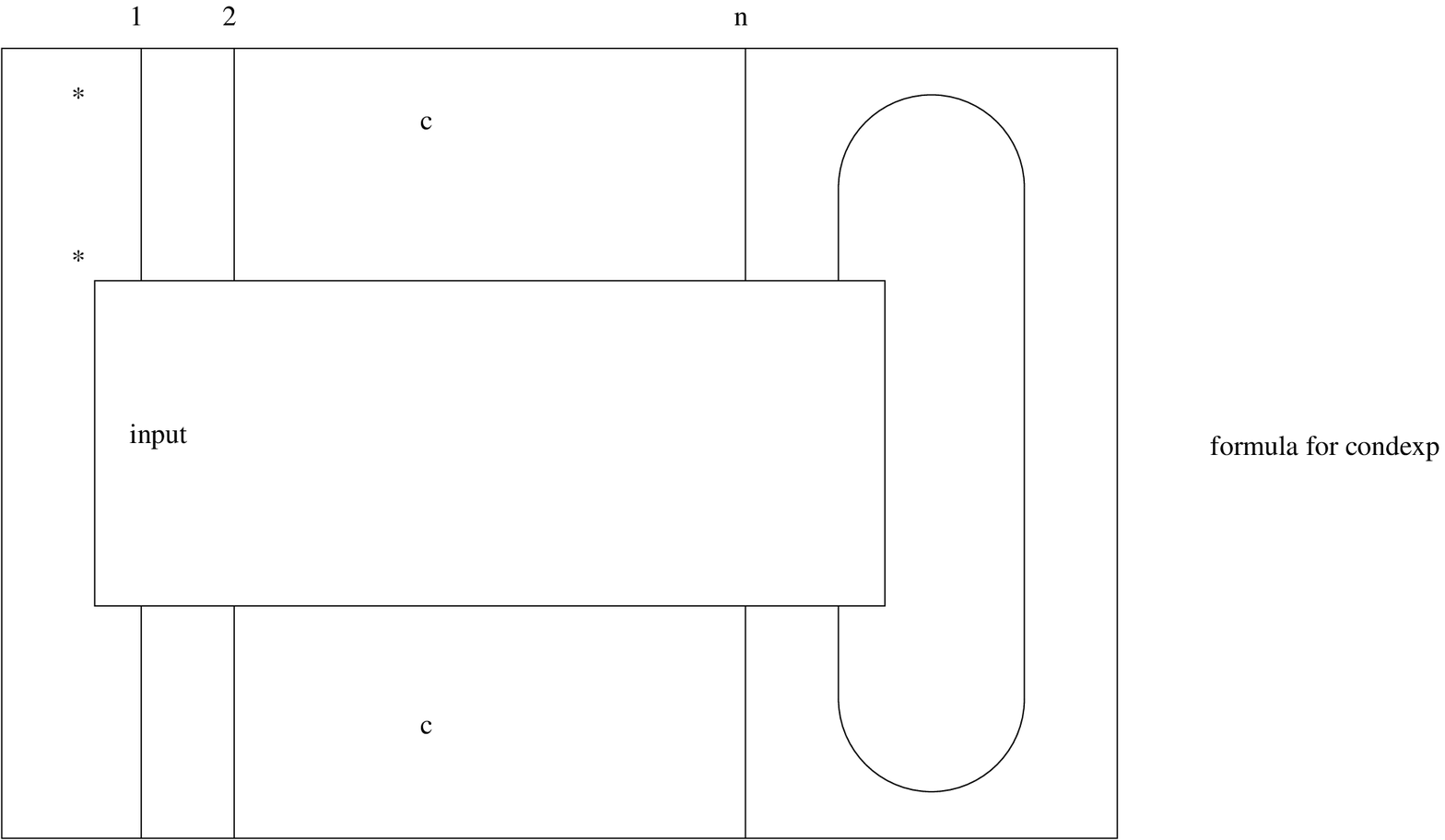, scale=0.3}
\end{flushleft}
and for $n = 0$, 
\psfrag{formula for 1stcondexp}
{$= \delta_{h,e} \delta_{k,e} \sqrt{\frac{\abs{K}}{\abs{H}}} e$}
\begin{flushleft}  
      \epsfig{file=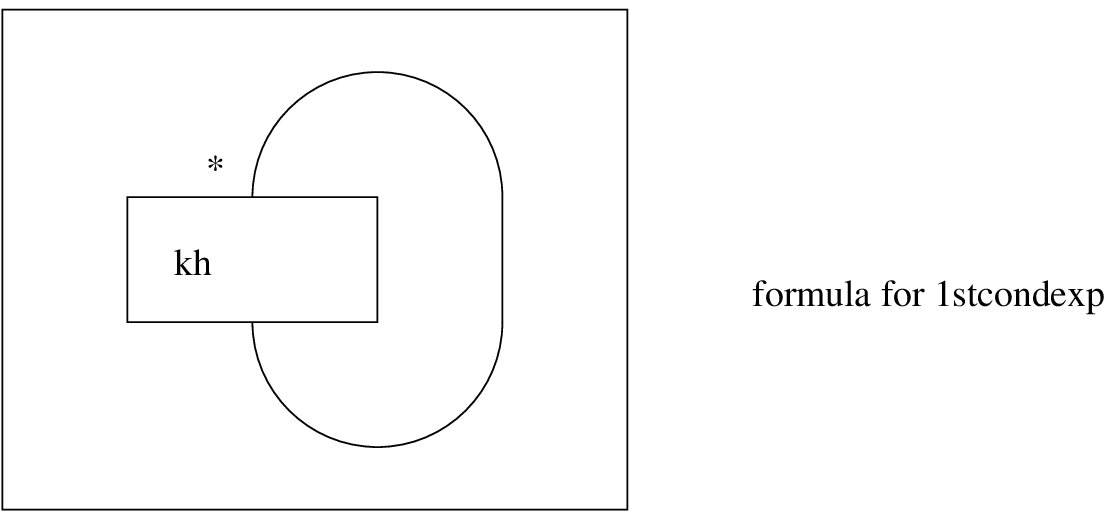, scale=0.5}
    \end{flushleft}
\item[(vii)] {\em Conditional Expectation Tangle from $P_n$ onto
  $P_{1,n}$}:\\ 
For $n \geq 2$ let $k_1 , k_2 \in K$, $\ul{t} \in T_{2n-3}$, 
$h \in H$ such that $\mu(k_1 , \ul{t} , k_2 , h) = e$. Then, 
\psfrag{leftinput}
{$(k_1 , \ul{t} , k_2 , h)$}
\psfrag{formula for leftcondexp}
{$ = \delta_{h,e} 
\sqrt{\frac{\abs{H}}{\abs{K}}}
\sum_{\substack{k', k'' \in K\\ \text{s.t.} k'' k' = k_2 k_1}}
(k', \ul{t}, k'', e)$}
\begin{flushleft}  
      \epsfig{file=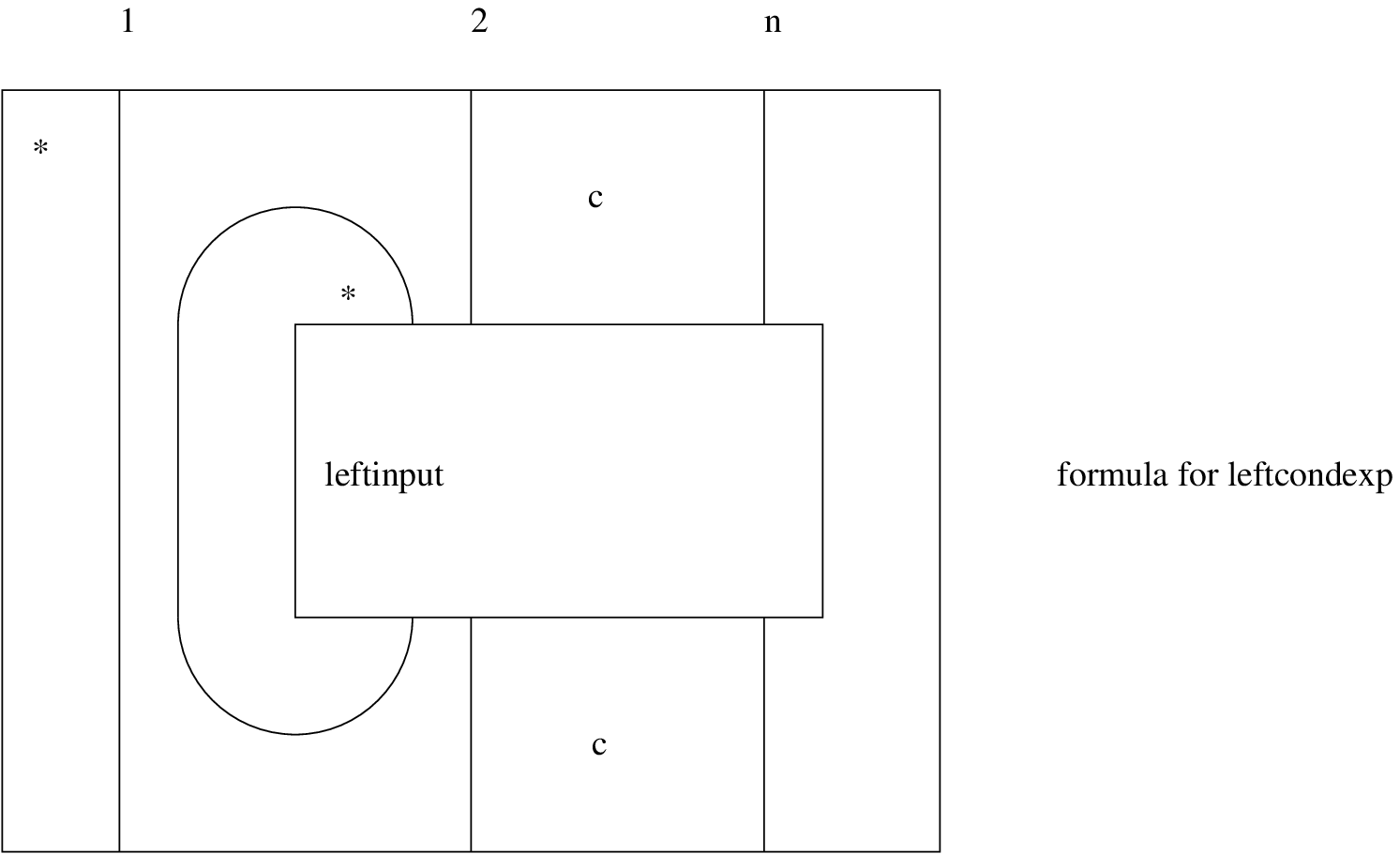, scale=0.4}
    \end{flushleft}
and for $n=1$, 
\psfrag{formula for 1stleftcondexp}
{$= \delta_{h,e} \delta_{k,e} \sqrt{\frac{\abs{H}}{\abs{K}}} e$}
\begin{flushleft}  
      \epsfig{file=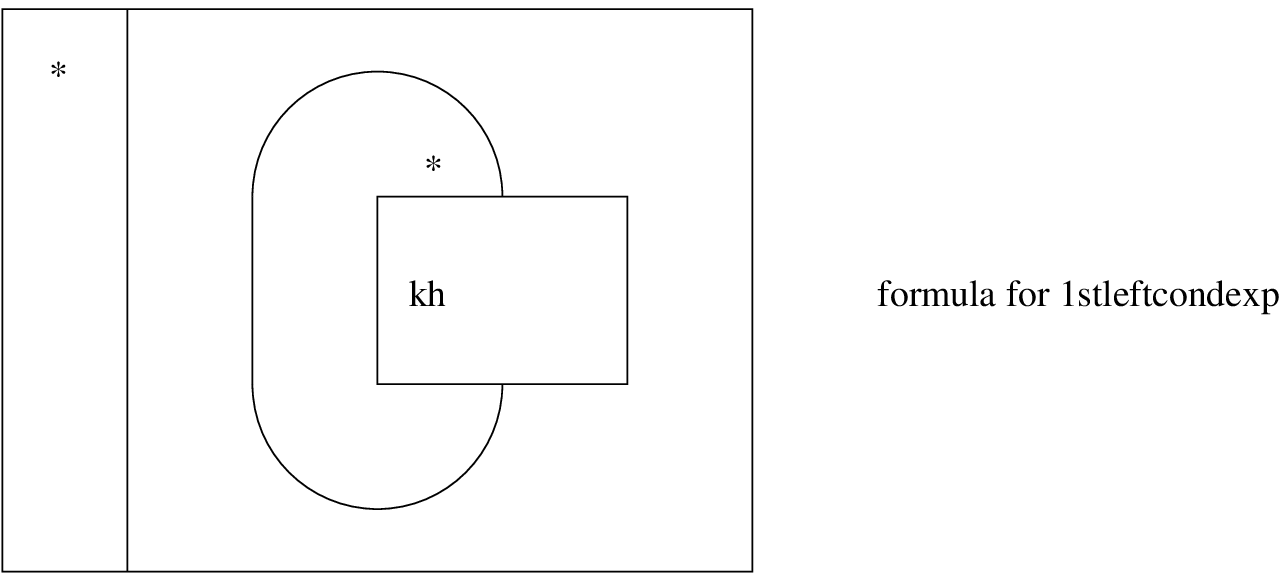, scale=0.4}
    \end{flushleft}
\end{itemize}
\end{remark}
\section{Iterated basic construction and higher relative commutants
of group-type subfactors}

Let $P$ be a $II_1$ factor and $G$ a countable discrete group
with an outer action on $P$, and suppose
$G$ is generated by two of its finite
subgroups $H$ and $K$. Consider the associated group-type subfactor
$P^H \subset P \rtimes K$.  
In this section we will give a concrete realization of the Jones tower 
associated to this subfactor and compute its higher 
relative commutants. See also \cite{Va} for related results.

First, let us recall the following characterization 
of the basic construction of a finite index subfactor (\cite{PiPo}).
\begin{lemma}\label{char}
Let $N \subset M$ be a finite index subfactor with $\E : M \rightarrow
N$ being the trace-preserving conditional expectation, $B$ be a $II_1$ factor
containing $M$ as a subfactor and $f$ be a projection in $B$ satisfying:
\begin{itemize}
\item[(i)] $fxf = \E(x)f$ for all $x \in M$,
\item[(ii)] $B$ is the algebra generated by $M$ and $f$.
\end{itemize}
Then, $B$ is isomorphic to the basic construction $M_1$ of $N \subset M$. 
\end{lemma}

Second, we recall some basic facts and
notations for the crossed product construction. 
Unless otherwise specified, we will reserve the symbol $e$ for 
the identity element of a group. The crossed product $P \rtimes K$
can be realized as the von Neumann subalgebra of
$\mathcal L(l^2(K) \otimes L^2(P))$
($\cong M_K (\mathbb{C}) \otimes \mathcal L(L^2(P))$) generated by the images
of $P$ and $K$ in the following way:
\begin{equation}\label{Pincl}
P \ni x \mapsto \sum_{k \in K} E_{k,k} \otimes k^{-1}(x) \in M_K (\mathbb{C}) \otimes \mathcal L(L^2(P))
\end{equation}
\begin{equation}\label{Kincl}
K \ni k \mapsto \lambda_k \otimes 1 \in M_K (\mathbb{C}) \otimes \mathcal L(L^2(P))
\end{equation}
where we set the convention of considering $k(x)$ as the element of $P$
obtained by applying the automorphism $k$ on $x$ (in $P$) and $\lambda_k$
is the matrix in $M_K (\C)$ corresponding to left multiplication by
$k$. Consequently, the following commutation relation holds in $P \rtimes K$:
\begin{equation}\label{algstr}
kxk^{-1} = k(x) \text{ for all } x \in P, k \in K.
\end{equation}
However, $P \rtimes K$ can also be viewed as the vector space generated by
elements of the form $\sum_{k \in K} x_k k$ where $x_k \in P$ where the
multiplication structure is given by the relation (\ref{algstr}).
The unique trace on $P \rtimes K$ is given by:
\[ tr(\sum_{k \in K} x_k k) = tr(x_e)\]  
and the unique trace-preserving conditional expectation is given by:  
\[\E^{P \rtimes K}_P (\sum_{k \in K} x_k k) = x_e.\]
If $P^K$ denotes the fixed point subalgebra of $P$, then $P \rtimes K$ is
isomorphic to the basic construction of $P^K \subset P$ where the Jones
projection is given by: 
\[ e_1 = \frac{1}{\abs K} \sum_{k \in K} k ~~ \in P\rtimes K\]
implementing the conditional expectation:
\[ \E^P_{P^K} (x) = \frac{1}{\abs K}\sum_{k \in K} k(x) ~~ \in P^K
\text{ for all } x \in P.\]
The basic construction of $P \subset P \rtimes K$ is isomorphic to
$M_K(\C) \otimes P$
where the inclusion $P \rtimes K \hookrightarrow M_K(\C)
\otimes P$ is given by the maps (\ref{Pincl}), (\ref{Kincl}),
and the corresponding Jones projection is given by 
$e_2 = E_{e,e} \otimes 1  \in M_K(\C) \otimes P$
implementing the conditional expectation 
$\E^{P \rtimes K}_P$.
The next element in the tower of basic construction is given by 
$M_K(\C) \otimes (P \rtimes K)$ where the inclusion
$M_K(\C) \otimes P \hookrightarrow M_K(\C) \otimes (P \rtimes K)$ is induced
by the inclusion $P \subset P \rtimes K$ and the Jones projection is given by:
\begin{equation}\label{e3}
e_3 = \frac{1}{\abs K} \sum_{k \in K} \rho_{k^{-1}} 
\otimes k ~~ \in M_K(\C) \otimes (P \rtimes K)
\end{equation}
implementing the conditional expectation:
\begin{equation}\label{condexp3}
\E^{M_K(\C) \otimes P}_{P \rtimes K} (E_{k_1,k_2} \otimes x) =
\frac{1}{\abs K} k_1 x {k_2}^{-1} ~~ \in P \rtimes K \text{ for all }
x \in P, k_1, k_2 \in K
\end{equation}
where $\rho_k$ is the matrix in $M_K (\C)$ corresponding to
right multiplication by $k$.

Coming back to the context of group-type subfactors 
we consider the unital inclusions $ P^H \hookrightarrow 
P\rtimes K  \hookrightarrow M_K(\C) \otimes (P \rtimes H)$ where the
second inclusion factors through $M_K(\C) \otimes P$ in the obvious way.
\begin{lemma}\label{firstlevel}
$M_K(\C) \otimes (P \rtimes H)$ is the basic construction for  
$P^H \subset P \rtimes K$ with Jones projection 
$e_1 = E_{e,e} \otimes \frac{1}{\abs H} \sum_{h \in H} h$.
\end{lemma}
\noindent {\bf Proof}:
We need to show that conditions (i) and (ii) of Lemma \ref{char} 
are satisfied.
To show (i), let us assume that $\tilde{x} = \sum_{k \in K} x_k k$
denotes a typical element of $P \rtimes K$.
\begin{eqnarray*}
e_1 \tilde{x} e_1 & = &
\Big(E_{e,e} \otimes \frac{1}{\abs H} \sum_{h \in H} h \Big) \; \tilde{x}
\; \Big(E_{e,e} \otimes \frac{1}{\abs H} \sum_{h{'} \in H} h{'} \Big) \\
& = & \Big(E_{e,e} \otimes \frac{1}{\abs H} \sum_{h \in H} h \Big) \; 
\Big(\sum_{\substack{k{'} \in K \\k \in K}}E_{k{'},k{'}}\lambda_k
\otimes {k{'}}^{-1}(x_k) \Big)
\; \Big(E_{e,e} \otimes \frac{1}{\abs H} \sum_{h{'} \in H} h{'} \Big)\\
& = & 
\Big(\sum_{k \in K}E_{e,e} \lambda_k E_{e,e} \otimes 
\frac{1}{{\abs H}^2} \sum_{\substack{h \in H \\h{'} \in H}}h x_k h{'} \Big)\\
& = & 
\Big(E_{e,e} \otimes \frac{1}{{\abs H}^2} 
\sum_{\substack{h \in H \\h{'} \in H}} h(x_e)hh{'} \Big)
\end{eqnarray*}
whereas
\begin{eqnarray*}
\E(\tilde{x}) e_1 & = &
\Big(\sum_{k \in K}E_{k,k} \otimes 
\frac{1}{\abs H} \sum_{h \in H} k^{-1}(h(x_e)) \Big)
\Big(E_{e,e} \otimes \frac{1}{\abs H} \sum_{h{'} \in H} h{'} \Big)\\
&=& \Big(E_{e,e} \otimes \frac{1}{{\abs H}^2} 
\sum_{\substack{h \in H\\h{'} \in H}} h(x_e) h{'} \Big)
\end{eqnarray*}
Therefore, LHS = RHS. 

To show (ii), it is enough to show that elements of the form 
$E_{k_1, k_2} \otimes xh$ for $x \in P,~ h \in H,~ k_1, k_2 \in K$
are in the algebraic span of $P \rtimes K$ and $e_1$. Let us denote 
the Jones projection in $P \rtimes H$ corresponding to the
inclusion $P^H \subset P$ by $f = \frac{1}{\abs H} \sum_{h \in H}
h$. Thus, \[e_1 = E_{e,e} \otimes f ~~ \in M_K(\C) \otimes P \rtimes H.\]
This implies \[P e_1 P = E_{e,e} \otimes PfP = E_{e,e} \otimes P
\rtimes H \subset M_K(\C) \otimes P \rtimes H\] where $P$ in the
left hand side is identified with its image inside  
$M_K(\C) \otimes P \rtimes H$ (namely, the prescription given by
(\ref{Pincl})).
Thus the algebraic span of $P \rtimes K$ and $e_1$ contains elements 
of the type $E_{e,e} \otimes xh$ for $x \in P$, $h \in H$. To obtain 
elements of the form $E_{k_1,k_2} \otimes xh$, note that the relation 
$\lambda_{k_1}E_{e,e}\lambda_{k_2^{-1}} = E_{k_1, k_2}$ holds 
in $M_K(\C)$.\qed
\begin{lemma}\label{secondlevel}
$M_K(\C) \otimes M_H(\C) \otimes (P \rtimes K)$ is the basic construction for  
$P \rtimes K \subset M_K(\C) \otimes (P \rtimes H)$
where the Jones projection is given by 
$e_2 = \frac{1}{\abs K} \sum_{k \in K} \rho_{k^{-1}} \otimes E_{e,e}  
\otimes k$.
\end{lemma}
\noindent {\bf Proof}: 
The conditional expectation 
$\E^{M_K(\C) \otimes (P \rtimes H)}_{P \rtimes K}$ is the composition 
$\E^{M_K(\C) \otimes P}_{P \rtimes K} \circ
\E^{M_K(\C) \otimes (P \rtimes H)}_{M_K(\C) \otimes P}$.
Therefore, $\E(E_{k_1, k_2} \otimes xh) = \delta_{h,e}
\frac{1}{\abs K} k_1x{k_2}^{-1}$.

\medskip 
To show condition (i) of Lemma \ref{char},  
\begin{eqnarray*}
&& e_2(E_{k_1, k_2} \otimes xh)e_2 \\
&=& 
\frac{1}{{\abs K}^2}\Big(\sum_{k{'} \in K} \rho_{k{'}^{-1}}\otimes 
E_{e,e} \otimes k{'}\Big)
\Big( \sum_{h{'} \in H} E_{k_1, k_2}\otimes E_{h{'},h{'}}
\lambda_h \otimes h{'}^{-1}(x) \Big)\\
&\ \ \ &
\Big(\sum_{k{''} \in K}\rho_{k{''}^{-1}}\otimes E_{e,e}\otimes
k{''}\Big)\\
&=& \delta_{h,e} \frac{1}{{\abs K}^2}\Big(
\sum_{k{'},k{''} \in K} \rho_{k{'}^{-1}} E_{k_1, k_2}
\rho_{k{''}^{-1}} \otimes E_{e,e} \lambda_h E_{e,e} \otimes
k{'}xk{''}\Big)\\ 
&=& \delta_{h,e} \frac{1}{{\abs K}^2}\Big(
\sum_{k{'},k{''} \in K} E_{k_1 k{'}^{-1},k_2 k{''}}
\otimes E_{e,e} \otimes k{'}xk{''}\Big)
\end{eqnarray*}
whereas 
\begin{eqnarray*}
&& \E(E_{k_1, k_2} \otimes xh)e_2\\
&=& \delta_{h,e} \frac{1}{{\abs K}^2} \Big(\sum_{h{'} \in H} 
\sum_{k \in K}\lambda_{k_1} E_{k,k} \lambda_{k_2^{-1}} \otimes
E_{h{'}, h{'}} \otimes h{'}^{-1}(k^{-1}(x))
\Big)
\Big(\sum_{k{'} \in K} \rho_{k{'}^{-1}} \otimes E_{e,e} \otimes k{'}\Big)\\
&=& \delta_{h,e} \frac{1}{{\abs K}^2} \Big(
\sum_{k \in K}\lambda_{k_1} E_{k,k} \lambda_{k_2^{-1}} 
\rho_{k{'}^{-1}} \otimes E_{e,e} \otimes k^{-1}(x)k{'}\Big)\\
&=& \delta_{h,e} \frac{1}{{\abs K}^2} \Big(
\sum_{k, k{'} \in K} E_{k_1 k, k_2 k k{'}} \otimes E_{e,e} \otimes
k^{-1}xkk{'} \Big)
\end{eqnarray*}
and the two sides are the same after renaming the indices. 

To show condition (ii), note that 
$M_K(\C) \otimes (P \rtimes K)$ is algebraically generated by 
$M_K(\C) \otimes P$ and $\frac{1}{\abs K} \sum_{k \in K} \rho_{k^{-1}} 
\otimes k$ ~(by the remarks preceding Equation (\ref{e3})). 
The following holds in 
$M_K(\C) \otimes M_H(\C) \otimes (P \rtimes K)$ 
because of this fact and the way $P$ sits inside $M_H(\C) \otimes P$
(since in the second tensor component we get expressions of the form
$\sum_{h, h{'} \in H}E_{h,h} E_{e,e} E_{h{'}, h{'}}$ 
which reduces to $E_{e,e}$):
\begin{eqnarray*}
(M_K(\C) \otimes P)e_2(M_K(\C) \otimes P) &=& M_K(\C) \otimes E_{e,e}
\otimes (P \rtimes K)\\ \Rightarrow 
(M_K(\C) \otimes P \rtimes H)e_2(M_K(\C) \otimes P \rtimes H) 
&=& M_K(\C) \otimes M_H(\C) \otimes (P \rtimes K)
\end{eqnarray*}
where the last implication is again due to the relation 
$\lambda_{h_1}E_{e,e}\lambda_{h_2^{-1}} = E_{h_1, h_2}$ in $M_H(\C)$.
\qed
\vskip 1em
Thus, we have the first two levels in the tower of basic construction:
\[P^H \subset P \rtimes K \subset M_K(\C) \otimes (P \rtimes H)
\subset M_{K \times H}(\C) \otimes (P \rtimes K) \]
where we identify $M_{K \times H}(\C)$ with $M_K(\C) \otimes M_H(\C)$. 
The next levels in the tower are obvious generalizations and we gather
everything in the following proposition.
\begin{proposition}\label{tower}
Let $G$ be a group acting outerly on the $II_1$
factor $P$ and assume $G$ is generated by two of its finite subgroups
$H$ and $K$. Then the $n$-th element of the tower of basic construction of the
group-type subfactor $N= P^H \subset P \rtimes K = M$ is given by:
\[M_n \cong M_{S_n}(\C) \otimes (P \rtimes L_n) \]
where the inclusion of $M_n$ inside $M_{n+1}$ is as follows:
\[M_n \ni E_{\ul{s},\ul{t}}\otimes x \mapsto \sum_{l \in L_n}
E_{\ul{s},\ul{t}} \otimes E_{l,l} \otimes l^{-1}(x) \in M_{n+1}
\text{ for all } x \in P , \; \ul{s} , \ul{t} \in S_n \]
\[M_n \ni E_{\ul{s},\ul{t}}\otimes l \mapsto
E_{\ul{s},\ul{t}} \otimes \lambda_l \otimes e \in M_{n+1}
\text{ for all } l \in L_n , \; \ul{s} , \ul{t} \in S_n \]
and the $n$-th Jones projection is:
\[
M_n \ni e_n = \left \{
\begin{array}{ll}
 \frac{1}{\abs {L_n}} \sum_{l \in L_n} I_{M_{S_{n-2}}} \otimes
\rho_{l^{-1}} \otimes E_{e,e} \otimes l, & \text{if } n > 1\\
\frac{1}{\abs H} \sum_{h \in H} E_{e,e} \otimes h, & \text{if } n = 1~ .
\end{array} 
\right.
\]
\end{proposition}
\noindent {\bf Proof}: We use induction. The case of $n=1$ is
a little different from the rest and is proved in Lemma
\ref{firstlevel} and the $n=2$ case is proved in Lemma
\ref{secondlevel}. Suppose the statement of the
above proposition holds upto a level $n$ ($> 2$). Now, the
subfactor $M_{n-1} \subset M_n$ is isomorphic to:
\[
M_{S_{n-1}} \otimes P \rtimes L_{n-1} \subset M_{S_{n-1}} \otimes
M_{L_{n-1}} \otimes P \rtimes L_n
\]
where we identify $M_{S_{n-1}} \otimes
M_{L_{n-1}}$ with $M_{S_n}$ and the inclusion is induced by identity
over $M_{S_{n-1}}$ tensored with the inclusion of the subfactor $P
\rtimes L_{n-1} \subset M_{L_{n-1}} \otimes P \rtimes L_n$. Using
Lemma \ref{secondlevel} for $K = L_{n-1}$ and $H = L_n$, it is clear
that the statement of the proposition holds for level $n+1$. \qed
\begin{remark}\label{condexprelcomm}
The formula for the unique trace-preserving 
conditional expectation is:
\[
\E^{M_n}_{M_{n-1}} (E_{\ul{s_1},\ul{s_2}} \otimes E_{m_1,m_2} \otimes
xl) = \delta_{l,e} \frac{1}{\abs {L_n}}
E_{\ul{s_1},\ul{s_2}} \otimes m_1 x m_2^{-1}
\]
where $\ul{s_1},\ul{s_2} \in S_{n-1}$, $m_1, m_2 \in L_{n-1}$, $l \in
L_n$ and $x \in P$ and the unique trace on $M_n$ is given by:
\[ tr_{M_n}(E_{\ul{s_1}, \ul{s_2}} \otimes x l)
= \frac{1}{\abs {S_n}} {\delta}_{l,e} {\delta}_{\ul{s_1}, \ul{s_2}} tr_M(x)
\]
where $\ul{s_1},\ul{s_2} \in S_n$, $l \in L_n$ and $x \in P$. 
\end{remark}
\vskip 1em
We will now compute the higher relative commutants using the above model 
of the Jones tower. To this end, we need the following two lemmas, where
we denote the set of automorphisms of $M \supset N$ that fixes
elements of $N$ pointwise by $Gal(N \subset M)$.
\begin{lemma}\label{galaut1}
Let $N \subset M$ be a finite index subfactor and $\theta \in 
Gal(N \subset M)$,  then the bimodule $_{M} L^2
(\theta) _M$ (where the module is $L^2(M)$ with usual left action of
$M$ but right action is twisted by $\theta$) is a $1$-dimensional
irreducible sub-bimodule of $_{M} L^2(M_1) _M$.
\end{lemma}
\noindent {\bf Proof}: Define $u_\theta : L^2(M) \rightarrow L^2(M)$
by $u_\theta (x\Omega) = \theta(x) \Omega$ for $x \in M$ where $\Omega$
is the cyclic and separating vector in the GNS construction with
respect to the canonical trace $tr$. Note that $u_\theta (n_1 \cdot
x\Omega \cdot n_2) = u_\theta (n_1 x n_2 \Omega) = n_1 \cdot \theta(x)
\Omega \cdot n_2$ for $n_1 , n_2 \in N$, $x \in M$. This implies
$u_\theta \in N' \cap
M_1$. Now define ${T :} _M L^2 (\theta)_M ~ {\rightarrow} ~ _M L^2(M_1)_M$ by
$T(x \Omega) = x u_\theta \Omega_1$ for $x \in M$. It is completely
routine to check that $T$ is a well-defined $M$-$M$ linear isometry and
we leave this to the reader. \qed
\begin{corollary}\label{galaut2}
$H = Gal(P^H \subset P)$ 
\end{corollary}
\noindent {\bf Proof}: Clearly $H \subset Gal(P^H \subset P)$. 
Let $\theta \in Gal(P^H \subset P)$. Note that $_P L^2(P \rtimes H)_P
\cong \oplus_{h \in H} ~_{P}L^2(h)_P$. Thus by Lemma \ref{galaut1},
$_P L^2(\theta)_P \, {\cong} \, _P L^2(h)_P$ for some $h \in H$. This
implies $\theta h^{-1} \in Inn(P) \cap Gal(P^H \subset P) = \{id_P\}$
since $P^H \subset P$ is irreducible. Hence, $\theta \in H$. \qed
\begin{lemma}\label{galaut3}
Let $N \subset M$ be an irreducible subfactor, i.e $N' \cap M \cong \C$
and $\theta \in Aut(M)$. For $x \in M$, the following are equivalent:
\begin{itemize}
\item[(i)] $x \neq 0$ and $x \theta(y) = y x$ for all $y \in N$,
\item[(ii)] $x_0 := \frac{x}{\Arrowvert x \Arrowvert} \in
   {\mathcal U}(M)$ and  $Ad_{x_0} \circ \theta \in Gal(N \subset M)$.
\end{itemize}
\end{lemma}
\noindent {\bf Proof} (ii)$\Rightarrow$(i) part is easy.

For (i)$\Rightarrow$(ii), note that we also have $\theta(y) x^* = x^* y$ for
all $y \in N$. Thus, $x^* x \in \theta(N)' \cap M$ and $x x^* \in N'
\cap M$. Since $N' \cap M \cong \C$, $x x^* = x^* x = \Arrowvert x
\Arrowvert^2$. Hence, $x_0 \in
   {\mathcal U}(M)$ and $x_0 ~ \theta(y) ~ {x_0}^* = y$ for all $y \in N$. This
   implies $Ad_{x_0}
   \circ \theta \in Gal(N \subset M)$. \qed
\begin{proposition}
For the group-type subfactor $N= P^H \subset P \rtimes K = M$, the
relative commutants $N' \cap M_n$ and $M' \cap M_n$ are given by:
\[ N' \cap M_n \cong 
\left\{
\begin{array}{ll}
\C, & \text{ if } n = -1\\
span
\left\{
E_{\ul{s_1}, \ul{s_2}} \otimes l 
\left|
\begin{array}{l}
\ul{s_1}, \ul{s_2} \in S_n,\\
l \in L_n, \\
\mu(\ul{s_1}) l \mu(\ul{s_2})^{-1} \in H 
\end{array}
\right.
\right\}, 
& \text{ if } n \geq 0
\end{array}
\right.\]
\[ M' \cap M_n \cong 
\left\{ 
\begin{array}{ll}
\C, & \text{ if } n = 0\\
span \left\{ 
\sum_{k \in K} E_{k, k k_0} \otimes E_{\ul{t_1}, \ul{t_2}}
\otimes l
\left| 
\begin{array}{l}
\ul{t_1}, \ul{t_2} \in T_{n-1},\\
\mu(\ul{t_1}) l \mu(\ul{t_2})^{-1} \in K,\\ 
k_0 = \mu(\ul{t_1}) l \mu(\ul{t_2})^{-1}
\end{array}
\right.
\right\},& \text{ if } n \geq 1.
\end{array}
\right.\]
\end{proposition}
\noindent {\bf Proof}: We compute the relative commutants in relation
to the concrete model of the
basic construction described in Proposition \ref{tower}. Consider the inclusion
\[N = P^H \ni x \mapsto \sum_{\ul{s} \in S_n}  E_{\ul{s}, \ul{s}}
\otimes \mu(\ul{s})^{-1} (x) \in M_{S_n}(\C) \otimes (P \rtimes L_n) =
M_n ~ .\]
\begin{eqnarray*}
\text {Let} && w = \sum_{\substack{\ul{s_1},\ul{s_2} \in S_n\\l \in L_n}}
E_{\ul{s_1}, \ul{s_2}} \otimes
x_{\ul{s_1},\ul{s_2}}^l l \in N' \cap M_n, \text{ for some } 
x_{\ul{s_1},\ul{s_2}}^l  \in P.\\
\text{ Then, } && w y  =  y w \text{ for all } y \in N\\
\Leftrightarrow && \sum_{\substack{\ul{s_1},\ul{s_2} \in S_n\\l \in L_n}}
E_{\ul{s_1}, \ul{s_2}} \otimes (x_{\ul{s_1},\ul{s_2}}^l \, l)
\mu(\ul{s_2})^{-1} (y) = \sum_{\substack{\ul{s_1},\ul{s_2} \in S_n\\l \in L_n}}
E_{\ul{s_1}, \ul{s_2}} \otimes \mu(\ul{s_1})^{-1} (y) \,
x_{\ul{s_1},\ul{s_2}}^l \, l \text{ for all } y \in N\\
\Leftrightarrow && x_{\ul{s_1},\ul{s_2}}^l (l \,
\mu(\ul{s_2})^{-1}) (y) = \mu(\ul{s_1})^{-1} (y) \,
x_{\ul{s_1},\ul{s_2}}^l \text{ for all } y \in N, \;
\ul{s_1},\ul{s_2} \in S_n, \; l \in L_n\\
\Leftrightarrow && \mu(\ul{s_1}) (x_{\ul{s_1},\ul{s_2}}^l) \, (\mu(\ul{s_1})
l \mu(\ul{s_2})^{-1}) (y) = y \, \mu(\ul{s_1}) (x_{\ul{s_1},\ul{s_2}}^l)
\text{ for all } y \in N, \; 
\ul{s_1},\ul{s_2} \in S_n, \; l \in L_n
\end{eqnarray*}
Now, by Lemma \ref{galaut3} and Corollary \ref{galaut2}, for
$\ul{s_1},\ul{s_2} \in S_n$ and $l \in L_n$, 
\[x_{\ul{s_1},\ul{s_2}}^l \neq 0 \Leftrightarrow
Ad_{x_0} \circ (\mu(\ul{s_1})
l \mu(\ul{s_2})^{-1}) \in H \text{ where } x_0 = \frac{\mu(\ul{s_1})
  (x_{\ul{s_1},\ul{s_2}}^l)}{\Arrowvert 
    \mu(\ul{s_1}) (x_{\ul{s_1},\ul{s_2}}^l) \Arrowvert} \; .\]
Moreover, $x_{\ul{s_1},\ul{s_2}}^l \neq 0 \Rightarrow Ad_{x_0} \in G
\cap Inn(P) = \{id_P\}$. Since $P$ is a factor,
$x_0 \in \C \, 1$. Thus, $x_{\ul{s_1},\ul{s_2}}^l$ will be
nonzero only if $\mu(\ul{s_1}) l \mu(\ul{s_2})^{-1} \in H$ and in such
cases $x_{\ul{s_1},\ul{s_2}}^l$ will be a scalar multiple of
identity. Hence, $N' \cap M_n$ is spanned by the linearly independent
set $\{E_{\ul{s_1}, \ul{s_2}} \otimes l : \ul{s_1}, \ul{s_2} \in S_n,
~ l \in L_n, ~ \mu(\ul{s_1}) l \mu(\ul{s_2})^{-1} \in H \}$.

For $M' \cap M_n$ where $n \geq 1$, we consider the inclusion
\begin{eqnarray*}
M = P \rtimes K \supset P \ni x & \mapsto & \sum_{\ul{s} \in S_n}
E_{\ul{s}, \ul{s}}\otimes \mu(\ul{s})^{-1} (x) \in M_{S_n}(\C) 
\otimes (P \rtimes L_n) = M_n \\
M = P \rtimes K \supset K \ni k & \mapsto & \lambda_k \otimes I_{T_{n-1}}
\otimes 1 \in M_{S_n}(\C) \otimes (P \rtimes L_n) = M_n 
\end{eqnarray*}
Let \[w = \sum_{\substack{\ul{s_1},\ul{s_2} \in S_n\\l \in L_n}}
E_{\ul{s_1}, \ul{s_2}} \otimes x_{\ul{s_1},\ul{s_2}}^l l \in M' \cap
M_n, \text{ for some } x_{\ul{s_1},\ul{s_2}}^l  \in P.\]

Using arguments similar to those for the calculations for the
case of $N' \cap M_n$ it follows that 
\begin{eqnarray*}
&& w y = y w \text{ for all } y \in P\\
\Leftrightarrow && \mu(\ul{s_1}) (x_{\ul{s_1},\ul{s_2}}^l) \, (\mu(\ul{s_1})
l \mu(\ul{s_2})^{-1}) (y) = y \, \mu(\ul{s_1}) (x_{\ul{s_1},\ul{s_2}}^l)
\text{ for all } y \in P, \; 
\ul{s_1},\ul{s_2} \in S_n, \; l \in L_n
\end{eqnarray*}
and hence $x_{\ul{s_1},\ul{s_2}}^l \neq 0 \; \Leftrightarrow \; 
\mu(\ul{s_1}) l \mu(\ul{s_2})^{-1} = e$ since $P$ is a factor; 
further, in such cases, $x_{\ul{s_1},\ul{s_2}}^l$ is a scalar multiple
of $1$. Now,
\begin{eqnarray*}
&& w k = k w \text{ for all } k \in K\\
\Leftrightarrow && k^{-1} w k = w \text{ for all } k \in K\\
\Leftrightarrow && \sum_{\substack{\ul{s_1},\ul{s_2}
    \in S_n,\\l \in L_n}} 
(\lambda_k^{-1} \otimes I_{T_{n-1}}
\otimes 1) (E_{\ul{s_1}, \ul{s_2}} \otimes
x_{\ul{s_1},\ul{s_2}}^l l) (\lambda_k \otimes I_{T_{n-1}}
\otimes 1)\\
&& =
\sum_{\substack{\ul{s_1},\ul{s_2} \in S_n\\l \in L_n}}
E_{\ul{s_1}, \ul{s_2}} \otimes
x_{\ul{s_1},\ul{s_2}}^l l \text{ for all } k \in K\\
\Leftrightarrow && \sum_{\substack{k_1, k_2 \in K,\\ \ul{t_1},\ul{t_2}
    \in T_{n-1},\\l \in L_n}} \lambda_k^{-1}
E_{k_1,k_2} \lambda_k \otimes E_{\ul{t_1}, \ul{t_2}} \otimes
x_{(k_1, \ul{t_1}),(k_2, \ul{t_2})}^l l\\
&& =
\sum_{\substack{k_1, k_2 \in K,\\ \ul{t_1},\ul{t_2}
    \in T_{n-1},\\l \in L_n}} 
E_{k_1,k_2} \otimes E_{\ul{t_1}, \ul{t_2}} \otimes
x_{(k_1, \ul{t_1}),(k_2, \ul{t_2})}^l l \text{ for all } k \in K\\
\Leftrightarrow && \sum_{\substack{k_1, k_2 \in K,\\ \ul{t_1},\ul{t_2}
    \in T_{n-1},\\l \in L_n}}
E_{k^{-1} k_1, k^{-1} k_2} \otimes E_{\ul{t_1}, \ul{t_2}} \otimes
x_{(k_1, \ul{t_1}),(k_2, \ul{t_2})}^l l\\
&& =
\sum_{\substack{k_1, k_2 \in K,\\ \ul{t_1},\ul{t_2}
    \in T_{n-1},\\l \in L_n}} 
E_{k_1,k_2} \otimes E_{\ul{t_1}, \ul{t_2}} \otimes
x_{(k_1, \ul{t_1}),(k_2, \ul{t_2})}^l l \text{ for all } k \in K\\
\Leftrightarrow && x_{(k_1, \ul{t_1}),(k_2, \ul{t_2})}^l = x_{(k k_1,
  \ul{t_1}),(k k_2, \ul{t_2})}^l \text{ for all } k, k_1, k_2 \in K,
\ul{t_1}, \ul{t_2} \in T_{n-1}, l \in L_n
\end{eqnarray*}
Finally, combining the conditions that we get from considering
commutation of $w$ with elements of $P$ and $K$, we can express $w$ as
a linear combination of elements of the form:
\[
\sum_{k \in K} E_{k k_1, k k_2} \otimes E_{\ul{t_1}, \ul{t_2}} \otimes l
\]
where $k_1, k_2 \in K$, $\ul{t_1}, \ul{t_2} \in T_{n-1}$, $l \in L_n$
such that $k_1 \mu(\ul{t_1}) l \mu(\ul{t_2})^{-1} k_2^{-1} =
e$. Equivalently, $w$ can be realized as a linear combination of:
\[
\sum_{k \in K} E_{k k_0^{-1}, k} \otimes E_{\ul{t_1}, \ul{t_2}}
\otimes l
= \sum_{k \in K} E_{k, k k_0} \otimes E_{\ul{t_1}, \ul{t_2}} \otimes l
\]
where $\ul{t_1}, \ul{t_2} \in T_{n-1}$, $l \in L_n$
such that $\mu(\ul{t_1}) l \mu(\ul{t_2})^{-1} \in K$ and $k_0 =
\mu(\ul{t_1}) l \mu(\ul{t_2})^{-1}$. \qed
\begin{remark}
(i) The set
$\left\{
E_{\ul{s_1}, \ul{s_2}} \otimes l 
\left|
\begin{array}{l}
\ul{s_1}, \ul{s_2} \in S_n,\\
l \in L_n, \\
\mu(\ul{s_1}) l \mu(\ul{s_2})^{-1} \in H 
\end{array}
\right.
\right\}$
(resp. \linebreak
$\left\{ 
\sum_{k \in K} E_{k, k k_0} \otimes E_{\ul{t_1}, \ul{t_2}}
\otimes l
\left| 
\begin{array}{l}
\ul{t_1}, \ul{t_2} \in T_{n-1},\\
\mu(\ul{t_1}) l \mu(\ul{t_2})^{-1} \in K,\\ 
k_0 = \mu(\ul{t_1}) l \mu(\ul{t_2})^{-1}
\end{array}
\right.
\right\}$
) forms a basis of $N' \cap M_n$ (resp. $M' \cap M_n$).\\

(ii)
The unique trace-preserving conditional expectation from $N' \cap M_n$
onto $M' \cap M_n$ is given by:
\[
\E^{N' \cap M_n}_{M' \cap M_n} (E_{k_1 , k_2} \otimes E_{\ul{t_1} ,
  \ul{t_2}} \otimes l) = \delta_{k_1 \mu(\ul{t_1}) l , k_2 
  \mu(\ul{t_2})} \frac{1}{\abs{K}} \sum_{k \in K} E_{k , k {k_1}^{-1} k_2}
  \otimes E_{\ul{t_1} , \ul{t_2}} \otimes l
\]
where $k_1 , k_2 \in K$, $\ul{t_1} , \ul{t_2} \in T_{n-1}$, $l \in
L_n$ such that $k_1 \mu(\ul{t_1}) l \mu(\ul{t_2})^{-1} k_2^{-1} \in H$.
\end{remark}
\noindent
To see (ii), we need to check that for 
$\ul{t_1}', \ul{t_2}' \in T_{n-1}$ and 
${k_0}' = \mu(\ul{t_1}') l' \mu(\ul{t_2}')^{-1} \in K$,
\begin{eqnarray*}
&& tr \Big [(\sum_{k' \in K} E_{k', k' {k_0}'} \otimes 
E_{\ul{t_1}', \ul{t_2}'} \otimes l')
(E_{k_1 , k_2} \otimes E_{\ul{t_1} ,\ul{t_2}} \otimes l)\Big ]\\
& = & \delta_{l'l, e}\; \delta_{\ul{t_2}', \ul{t_1}}\;
\delta_{\ul{t_1}', \ul{t_2}}\; \delta_{k_2 {k_0}', k_1} 
\frac{1}{\abs{S_n}}\\
& = & \delta_{k_1 \mu(\ul{t_1}) l , k_2 \mu(\ul{t_2})}\; 
\delta_{l'l, e}\; \delta_{\ul{t_1}, \ul{t_2}'}\;
\delta_{\ul{t_2}, \ul{t_1}'}\; \delta_{e, {k_0}' {k_1}^{-1} k_2} 
\frac{1}{\abs{S_n}}\\
& = & tr \Big [(\sum_{k' \in K} E_{k', k' {k_0}'} \otimes 
E_{\ul{t_1}', \ul{t_2}'} \otimes l')
(\delta_{k_1 \mu(\ul{t_1}) l , k_2 
  \mu(\ul{t_2})} \frac{1}{\abs{K}} \sum_{k \in K} E_{k , k {k_1}^{-1} k_2}
  \otimes E_{\ul{t_1} , \ul{t_2}} \otimes l) \Big ]
\end{eqnarray*}
which is a routine computation using the second part of 
Remark \ref{condexprelcomm}.
\section{The planar algebra of group-type subfactors}\label{main}
In this section, our main goal is to show that the planar algebra
defined in section \ref{plnalg}, is isomorphic to the one arising from the
group-type subfactor $P^H \subset P \rtimes K$. Conversely, we will
show that any subfactor whose standard invariant is given by such 
planar algebra, is indeed of this type.

We will use the following well-known fact regarding isomorphism of two
planar algebras and Theorem 4.2.1 of
\cite{J2} describing the planar algebra arising from an extremal
finite index subfactor.

\noindent{\bf Fact}: Let $P^1$ and $P^2$ be two planar algebras. Then $P^1
\cong P^2$ if and only if there exist a vector space isomorphism
$\psi : P^1 \rightarrow P^2$ such that:
\begin{itemize}
\item[(i)] $\psi$ preserves the filtered algebra structure,
\item[(ii)] $\psi$ preserves the actions of all possible Jones
  projection tangles and the (two types of) conditional expectation tangles.
\end{itemize}
If $P_1$ and $P_2$ are $*$-planar algebras, then we consider such $\psi$ that
are $*$-preserving to be a $*$-planar algebra isomorphism.

Let us denote the planar algebra in section \ref{plnalg} by $P^{BH}$
and the one arising from the group-type subfactor $N = P^H \subset P
\rtimes K = M$ by $P^{N \subset M}$.
\begin{theorem} \label{forward}
$P^{N \subset M} \cong P^{BH}$.
\end{theorem}
\noindent{\bf Proof}: By Theorem 4.2.1 of \cite{J2} we have that
$P^{N \subset M}_n = N' \cap M_n$ where the $n$-th Jones projection
tangle is given by $\delta e_n$, the conditional expectation tangle from
$P^{N \subset M}_n$ onto $P^{N \subset M}_{n-1}$ is $\delta \E^{N'
  \cap M_{n+1}}_{N' \cap M_n}$ and the conditional expectation tangle from
$P^{N \subset M}_n$ onto $P^{N \subset M}_{1, n}$ is $\delta \E^{N'
  \cap M_{n+1}}_{M' \cap M_{n+1}}$.
\vskip 1em
Define the map
\[
\begin{array}{rccc}
\psi : &  P^{N \subset M} & \rightarrow & P^{BH}\\
&&&\\
& \cup & & \cup\\
&&&\\
\psi_n : & P^{N \subset M}_n &\rightarrow & P^{BH}_n
\end{array}
\text{ by }
\psi_n (E_{\ul{s_1}, \ul{s_2}} \otimes l) = (\ul{s_1}, l,
\tilde{\ul{s_2}}, h) \]
where $n \geq 0$, $\ul{s_1} , \ul{s_2} \in S_n$, $l \in L_n$ such that
$\mu(\ul{s_1}) l \mu(\ul{s_2})^{-1} \in H$ and $\mu(\ul{s_1}) l
\mu(\ul{s_2})^{-1} h = e$.

Clearly, $\psi$ is a vector space isomorphism by definition. In order
to check that $\psi$ is a filtered $*$-algebra isomorphism, we use
Remark \ref{structure} (i), (ii), (iii) and (iv). For instance, to show
that $\psi_n$ is an algebra homomorphism, we need to show
\[
\begin{array}{ll}
&\psi_n (E_{\ul{s_1}, \ul{s_2}} \otimes l_1 \cdot E_{\ul{s_3},
  \ul{s_4}} \otimes l_2) = \psi_n (E_{\ul{s_1}, \ul{s_2}} \otimes l_1)
\; \cdot \; \psi_n (E_{\ul{s_3}, \ul{s_4}} \otimes l_2)\\
&\\
\Leftarrow & \delta_{\ul{s_2} , \ul{s_3}} \; (\ul{s_1}, l_1 l_2,
  \tilde{\ul{s_4}} , \mu(\ul{s_4}) l^{-1}_2 l^{-1}_1
  \mu(\ul{s_1})^{-1})\\
&\\
&=
(\ul{s_1}, l_1, \tilde{\ul{s_2}} , \mu(\ul{s_2}) l^{-1}_1
  \mu(\ul{s_1})^{-1}) \cdot
(\ul{s_3}, l_2, \tilde{\ul{s_4}} , \mu(\ul{s_4}) l^{-1}_2 \mu(\ul{s_3})^{-1})
\end{array}
\]
which indeed holds by Remark \ref{structure} (iii).

Now, it remains to be shown that $\psi$ preserves the action of Jones
projection tangles and the two types of conditional expectation
tangles. For this, we use Remark \ref{structure} (v), (vi) and
(vii). Proof of each of the three kinds of tangles, is completely
routine; however, we will discuss the action of conditional expectation
tangle in details. Let us consider the conditional expectation tangle
from $n+1$ (colour of the internal disc) to colour $n$ (colour of the
external disc) applied to the element $E_{\ul{s_1}, \ul{s_2}} \otimes
E_{m_1 , m_2} \otimes
l \in P^{N \subset M}_{n+1} = N^\prime \cap M_n$; by Jones's theorem
(4.2.1 of \cite{J2}) and Remark \ref{condexprelcomm}, the output
should be
\[
\begin{array}{ll}
& \delta_{l , e} \; \sqrt{\frac{\abs{L_n}}{\abs{L_{n+1}}}} \; E_{\ul{s_1},
    \ul{s_2}} \otimes m_1 m^{-1}_2\\
& \\
\stackrel{\psi_n}{\longmapsto} &
\delta_{l , e} \; \sqrt{\frac{\abs{L_n}}{\abs{L_{n+1}}}} \; (\ul{s_1},
m_1 m^{-1}_2, \tilde{\ul{s_2}}, \mu(\ul{s_2}) m_2 m^{-1}_1
    \mu(\ul{s_1})^{-1})
\end{array}
\]
On the other hand, by Remark \ref{structure} (vi), the conditional
expectation tangle applied to $\psi_{n+1} (E_{\ul{s_1}, \ul{s_2}} \otimes
E_{m_1 , m_2} \otimes l) = (\ul{s_1}, m_1, l^{-1}, m^{-1}_2,
\tilde{\ul{s_2}}, \mu(\ul{s_2}) m_2 l m^{-1}_1
    \mu(\ul{s_1})^{-1})$ is given by $\delta_{l , e} \; \sqrt{\frac{\abs{L_n}}{\abs{L_{n+1}}}} \; (\ul{s_1},
m_1 m^{-1}_2, \tilde{\ul{s_2}}, \mu(\ul{s_2}) m_2 m^{-1}_1
    \mu(\ul{s_1})^{-1})$.

This completes the proof.\qed
\begin{corollary}
Given any countable discrete group $G$ generated by two of its finite 
subgroups, there
exists a hyperfinite subfactor with standard invariant described by
$P^{BH}$. Moreover, $P^{BH}$ is a spherical $C^*$ planar algebra. 
\end{corollary}
\noindent{\bf Proof}: The proof follows from the fact that any countable
discrete group $G$ has an outer action on the hyperfinite $II_1$
factor.\qed
\begin{theorem}\label{converse}
Given a subfactor $N \subset M$ with standard invariant isomorphic to 
$P^{BH}$, there exists an intermediate subfactor $N \subset P \subset
M$ and outer actions of $H$ and $K$ on $P$ such that ($N \subset M$)
$\simeq$ ($P^H \subset P \rtimes K$).
\end{theorem}
\noindent{\bf Proof}: Let $P^{N \subset M}$ denote the planar algebra
of $N \subset M$ formed by its relative commutants and $\phi : P^{N
  \subset M} \rightarrow P^{BH}$ be a planar algebra
isomorphism. Consider the element $q = {\phi}^{-1}(e,e,e,e) \in N^\prime \cap
M_1$. Clearly, (i) $q$ is a projection, (ii) $q e_1 = e_1$ and (iii)
$\E_M (q) = \abs{K}^{-1} 1$. Using action of tangles and the planar
algebra isomorphism $\phi$, it also follows that
\psfrag{q}{$q$}
\psfrag{(iv)}{(iv)}\psfrag{1}{$1$}\psfrag{n}{$n$}
\psfrag{eqsqrthbyk}{$=\sqrt{\frac{\abs{H}}{\abs{K}}}$}
\begin{figure}[h]
\begin{center}  
      \epsfig{file=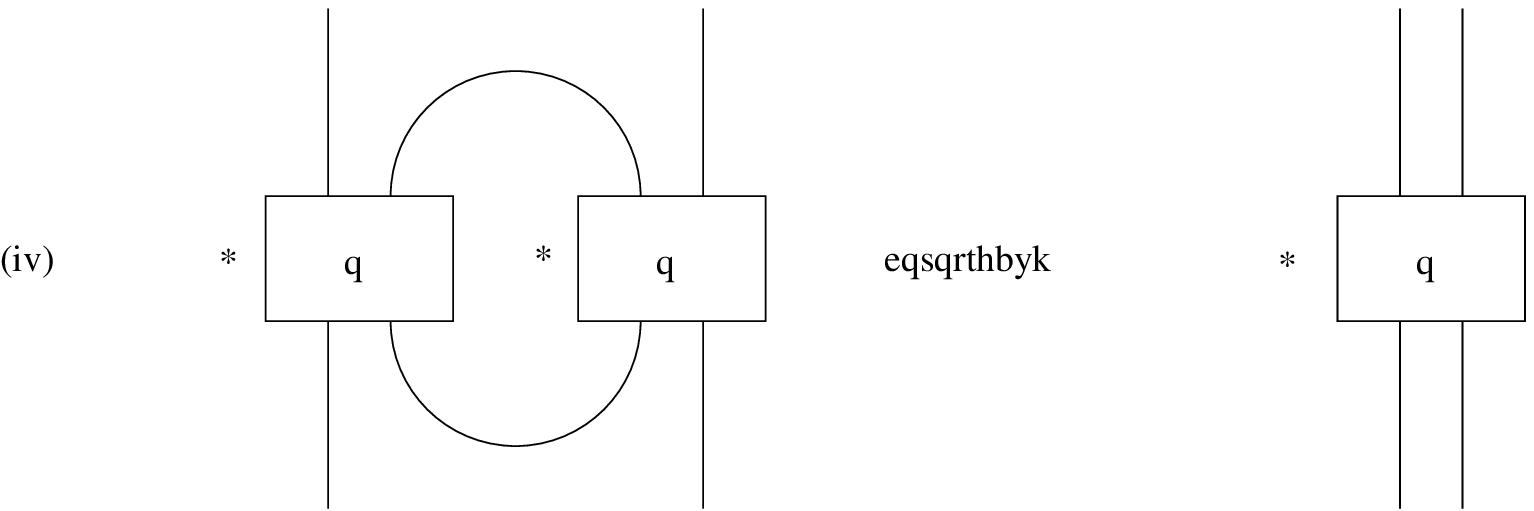, scale=0.4}
    \end{center}
\end{figure}\\
The conditions (i) - (iv) asserts that $q$ is an {\em intermediate
  subfactor projection} as described in \cite{Bi1}. Define
$P = M \cap \{q\}^\prime$. To show $P$ is a factor, first note that
\[P^\prime \cap P \subset N^\prime \cap P = N^\prime \cap M \cap
\{q\}^\prime = \phi(P^{BH}_1 \cap \{(e,e,e,e)\}^\prime)\]
So, it is enough to show that $P^{BH}_1 \cap \{(e,e,e,e)\}^\prime =
\C 1$. If $x = \sum_{g \in K \cap H} (g,g^{-1}) \in P^{BH}_1 \cap
\{(e,e,e,e)\}^\prime$, then
\[(e,e,e,e) \cdot x = x \cdot (e,e,e,e) \Rightarrow \sum_{g \in K \cap
  H} \lambda_g (e,e,g,g^{-1}) = \sum_{g \in K \cap H} \lambda_g
  (g,e,e,g^{-1})\]
Hence, $\lambda_g = \delta_{g,e} \lambda_e$ for all $g \in K \cap H$
  and $x \in \C 1$.

It remains to establish that $N$ (resp. $M$) is the fixed-point subalgebra
(resp. crossed-product algebra) of $P$ with respect to an outer action of
the group $H$ (resp. $K$). It is easy to prove (see \cite{J2}) that if
the standard
invariant of a subfactor $\tilde{N} \subset \tilde{M}$ is given by the
planar algebra corresponding to the fixed-point subfactor
(resp. crossed-product subfactor) with respect to a finite group
$\tilde{G}$, then there exist an outer action of $\tilde{G}$ on
$\tilde{M}$ (resp. $\tilde{N}$) such that $\tilde{N}$
(resp. $\tilde{M}$) is isomorphic to fixed-point subalgebra
(crossed-product algebra) of the action. Again, if $\tilde{N} \subset
\tilde{P} \subset \tilde{M}$ is an intermediate subfactor and
$\tilde{q}$ is its corresponding intermediate subfactor projection,
then the planar algebra of $\tilde{N} \subset \tilde{P}$
(resp. $\tilde{P} \subset \tilde{M}$) is given by the range of the
idempotent tangle
\psfrag{tq}{$\tilde{q}$}
\psfrag{c}{$\cdots$}
\psfrag{or}{or}
\begin{figure}[h]
\begin{center}  
      \epsfig{file=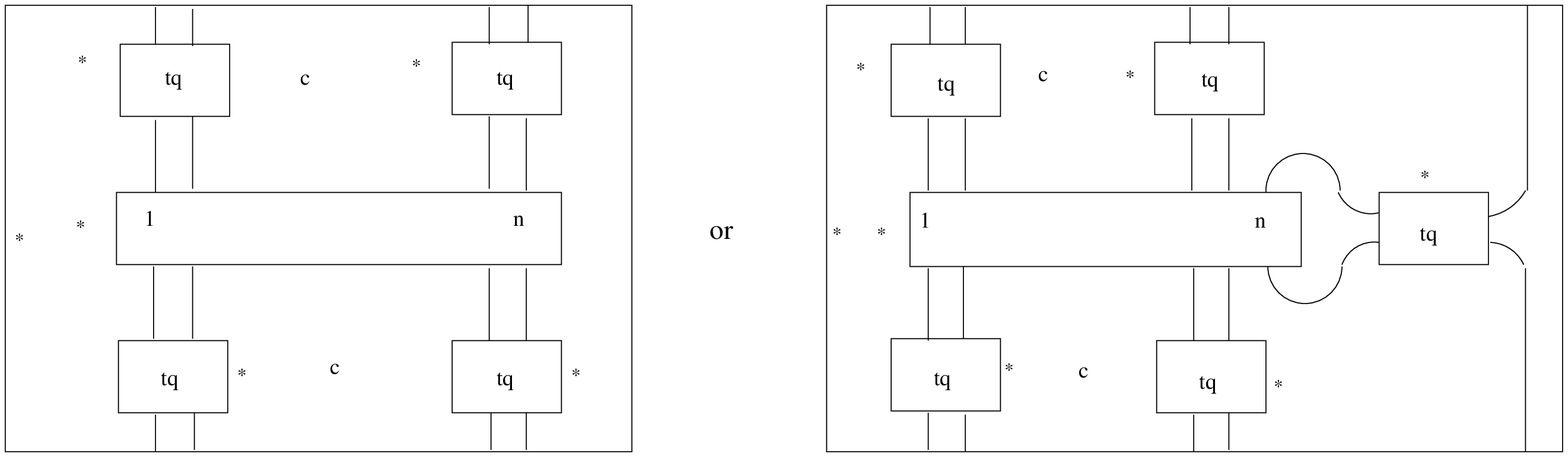, scale=0.4}
    \end{center}
\end{figure}
\psfrag{tq}{$\tilde{q}$}
\psfrag{c}{$\cdots$}
\psfrag{resp}{(resp.}
\psfrag{brcl}{)}
\begin{figure}[h]
\begin{center}  
      \epsfig{file=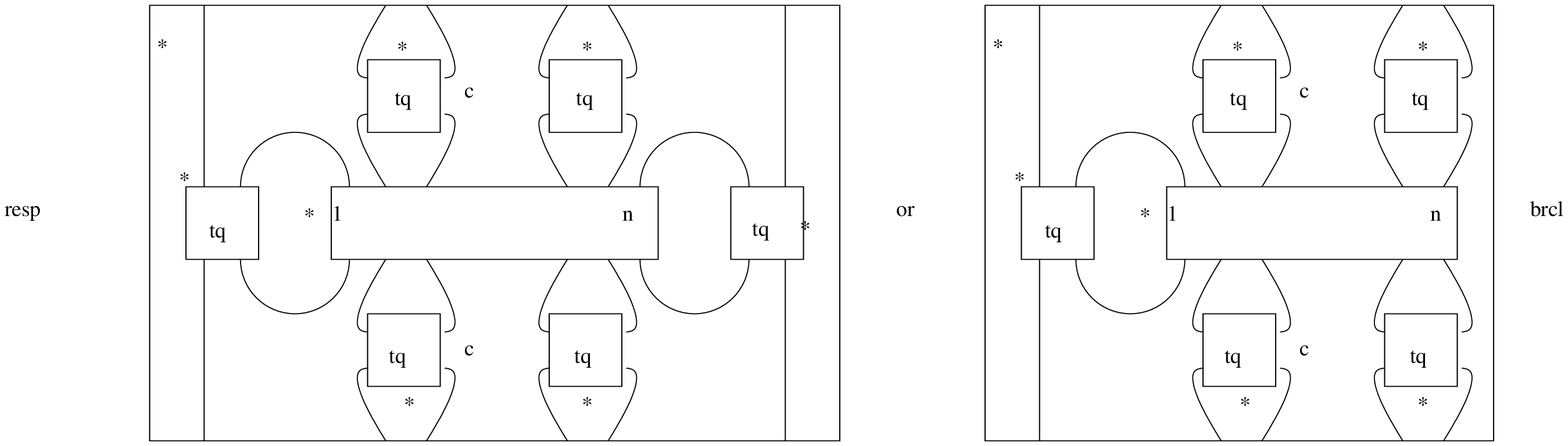, scale=0.45}
    \end{center}
\end{figure}\\
according as $n$ is even or odd (\cite{BJ2}, see also
\cite{BZ}).

Getting back to our context, to get the planar algebra of $N \subset
P$ (resp. $P
\subset M$), the elements in the image of the above idempotent tangles are
given by the words with letters coming from $K$
and $H$ alternately where every element coming from $K$ (respectively
$H$) must necessarily be $e$. Such a planar algebra is the same as
$P^{BH}$ with $K = \{e\}$ (resp. $H = \{e\}$); by Theorem
\ref{forward} this is indeed the planar algebra corresponding to
fixed-point subfactor (resp. crossed-product subfactor) with respect
to $H$ (resp. $K$).\qed
\bibliographystyle{amsplain}

\end{document}